\documentclass[12pt]{amsart}
\usepackage{amssymb}
\ifx\mathrm\undefined\let\mathrm\rm\fi
\ifx\mathbf\undefined\let\mathbf\bf\fi
\ifx\mathfrak\undefined\let\mathfrak\frak\fi
\ifx\mathcal\undefined\let\mathcal\cal\fi
\ifx\mathbb\undefined\let\mathbb\Bbb\fi
\ifx\emph\undefined\let\emph\it\fi
\newcommand{\g}{{{\mathfrak g}\,}}

\newcommand{\n}{{{\mathfrak n}}}
\newcommand{\h}{{{\mathfrak h\,}}}
\newcommand{\Tr}{{{\rm Tr}}}

\newcommand{\hr}{{{\mathfrak h_{\mathbb R}\,}}}

\newcommand{\Z}{{\mathbb Z}}

\newcommand{\R}{{\mathbb R}}
\newcommand{\C}{{\mathbb C}}

\newcommand{\Ref}[1]{{(\ref{#1})}}
\newcommand{\be}{\begin{displaymath}}
\newcommand{\ee}{\end{displaymath}}
\newcommand{\bea}{\begin{eqnarray*}}
\newcommand{\eea}{\end{eqnarray*}}

\newcommand{\bq}{\begin{equation}}

\newcommand{\eq}{\end{equation}}

\newcommand{\W}{{{\Bbb W\,}}}

\newcommand{\al}{{\alpha}}
\newcommand{\om}{{\omega}}
\newcommand{\dl}{{\delta}}
\newcommand{\Si}{{\Sigma}}

 at9.98pt

\newcommand{\dontprint}[1]{\relax}

\newcommand{\la}{\lambda}

\newcommand{\T}{\!\otimes\!\,}

\newcommand{\bean}{\begin{eqnarray}}
\newcommand{\eean}{\end{eqnarray}}

\newtheorem%
{thm}{Theorem}[section]
\newtheorem%
{proposition}[thm]{Proposition}
\newtheorem%
{lemma}[thm]{Lemma}

\newtheorem%
{lemmadef}[thm]{Lemma-Definition}
\newtheorem%
{corollary}[thm]{Corollary}
\newtheorem%
{conjecture}[thm]{Conjecture}

\title{
The orthogonality and $q$KZB-heat equation
for traces of $U_q(\g)$-intertwiners}
\author[
P. Etingof and A. Varchenko
]
{
Pavel Etingof${}^{*, 1}$  \and Alexander Varchenko${}^{**, 2}$
}

\begin{document}
\maketitle
\centerline{\it ${}^*$ Rm. 2-176, Department of Mathematics, MIT}
\centerline{\it Cambridge, MA, 02139, USA}
\centerline{etingof@math.mit.edu}
\medskip
\centerline{\it ${}^{**}$Department of Mathematics,
University of North Carolina at Chapel Hill,}
\centerline{\it Chapel Hill, NC 27599-3250, USA}
\centerline{anv@email.unc.edu}
\newcommand{\sig}{\sigma}
\newcommand{\Hom}{\text{Hom} \,}
\newcommand{\End}{\text{End} \,}
\newcommand{\IM}{\text{Im} \,}
\newcommand{\RE}{\text{Re} \,}

\vskip .1in
\centerline{\bf Dedicated to the 90-th birthday of Izrail Moiseevich Gelfand}

\begin{abstract} In our previous paper \cite{EV2},
to every finite dimensional representation $V$ of 
the quantum group $U_q({\frak g})$, we attached 
the trace function $F^V(\lambda,\mu)$,
with values in ${\rm End} \, V[0]$, obtained by taking the
(weighted) trace in a Verma module of an intertwining operator.
We showed that these trace functions satisfy the Macdonald-Ruijsenaars and the
qKZB equations, their dual versions, and the symmetry identity. 
In this paper we show that the trace functions satisfy the orthogonality
relation and the qKZB-heat equation. For ${\frak g}=sl_2$, 
this statement is the trigonometric degeneration 
of a conjecture from \cite{FV2}, proved in \cite{FV2} for
the 3-dimensional irreducible $V$. 

We also establish the orthogonality relation and qKZB-heat equation 
for trace functions obtained by taking traces in finite
dimensional representations (rather than Verma modules). 
If $\g=sl_n$ and $V=S^{kn}\Bbb C^n$, these 
functions are known to be Macdonald polynomials of type $A$. In this case, 
the orthogonality relation reduces to the Macdonald inner product identities, 
and the qKZB-heat equation coincides with the 
q-Macdonald-Mehta identity, proved by Cherednik \cite{Ch1}.  
\end{abstract}

\section{Introduction}
This paper is motivated by the previous papers
\cite{FTV1,FTV2,FV1,FV2,FV3,EV2}.

Let $V_1,...,V_n$ be finite dimensional representations of $U_q(sl_2)$. 
In the papers \cite{FTV1,FTV2}, G. Felder, V.Tarasov, and the second author introduced 
the function \linebreak $u_{V_1,...,V_n}(\lambda,\mu,\tau,p,q,z_1,...,z_n)$
with values in ${\rm End}((V_1\otimes...\otimes V_n)[0])$, 
which is a q-deformation of a conformal block
of the $sl_2$-Wess-Zumino-Witten conformal field theory 
on an elliptic curve with $n$ punctures. Namely, $u$ is defined
by a q-deformed version of the
explicit integral representation of conformal blocks 
on an elliptic curve. If $n=1$, and $V_1=L_{2m}$ is the representation with 
highest weight $2m$ (of dimension $2m+1$), then 
the function $u$ is independent of $z_1$ and scalar valued;
it is denoted by $u_m(\lambda,\mu,\tau,p,q)$. 

It is shown in \cite{FTV1,FTV2} that the function
$u$ satisfies the qKZB difference equations, their dual
version, and the symmetry with respect to the permutation 
$(\lambda,\tau)\leftrightarrow (\mu,p)$. 
Later G.Felder and the 
second author showed in \cite{FV2} that for $m=1$ the function $u_m$ 
satisfies the orthogonality relation
and the qKZB-heat equation. The latter is a q-deformation
of the usual Knizhnik-Zamolodchikov-Bernard (KZB) 
heat equation for conformal blocks. However, for $m>1$ the qKZB-heat equation 
still remains a conjecture, and the higher rank case is yet to be worked out. 

In order to understand the mysterious nature of the function
$u$ and in particular to answer the above questions, 
we proposed in \cite{EV2} to study the 
representation theoretic meaning of $u$. More specifically,
since conformal blocks on an elliptic
curve are known to be (weighted) traces of 
products of intertwining operators for affine
algebras, we conjectured that the function $u$ is obtained by taking 
a similar trace for the quantum affine algebra, and checked this
conjecture in the trigonometric limit.

More precisely, to every 
collection of finite dimensional representations $V_1,...,V_n$ of 
the quantum group $U_q({\frak g})$, we attached 
the trace function $F^{V_1,...,V_n}(\lambda,\mu)$, which depends on 
two complex weights $\lambda,\mu$ of $\g$ and the parameter $q$, and takes
values in ${\rm End}(V_1\otimes...\otimes V_n)[0]$, 
the endomorphism algebra of the zero
weight space of $V$. This function is by definition a suitably 
renormalized weighted 
trace in a Verma module of a product of intertwining operators.
The main results of \cite{EV2} are 
 that the trace functions satisfy the Macdonald-Ruijsenaars and the
qKZB equations, their dual versions, and the symmetry identity, 
and also that for ${\frak g}=sl_2$ and $V=L_{2m}$,
the function $F^V$ (up to simple renormalizations) coincides with the limit 
of $u_m$ as $p$ and $\tau$ go to infinity.  

In this paper we continue the study 
of the trace functions for $U_q(\g)$ (for $n=1$), 
and show that they satisfy the orthogonality
relation and the qKZB-heat equation. 
We also define the integral transform with 
kernel $F^V(\lambda,\mu)$. In the special case $V=\Bbb C$, 
this integral transform 
specializes to the usual Fourier transform, 
while for ${\frak g}=sl_n$, $V=S^{kn}\Bbb C^n$ 
its symmetrized version is the Cherednik's difference Fourier transform 
for type $A$ \cite{Ch2}. The orthogonality relation can be understood 
as the statement that the inverse to the 
integral transform with kernel $F^V(\mu,-\lambda)$ 
is the integral transform with the kernel
$F^V(\lambda,\mu)Q_V(-\mu-\rho)$
(where $Q$ is the contracted fusion operator, see \cite{EV2}), 
which generalizes the standard ``self-duality'' 
property of the Fourier and Cherednik transforms. 
For ${\frak g}=sl_2$, 
these results are the trigonometric degenerations of the statements 
conjectured in \cite{FV2}, and proved in  
the case when $V$ is 3-dimensional. 

We also establish the orthogonality relation and qKZB-heat equation 
for trace functions obtained by taking traces in finite
dimensional representations (rather than Verma modules). 
In the case $\g=sl_n$ and $V=S^{kn}\Bbb C^n$, these 
functions are Macdonald polynomials of type $A$ \cite{M}, 
the orthogonality relation reduces to the Macdonald inner product identities, 
and the qKZB-heat equation coincides with the 
q-Macdonald-Mehta identity, proved by Cherednik \cite{Ch1}.  

These results can be easily generalized to the case $n>1$, 
using the techniques of \cite{EV2}. We will not discuss this generalization. 

The structure of the paper is as follows. In Section 2 we 
recall the basics about quantum groups, the definition of trace functions 
from \cite{EV2}, and formulate our main results -- the orthogonality
and the q-KZB heat equation. We also 
reformulate the orthogonality 
relation as an inversion formula for a generalized Fourier
transform, similarly to how it was done in \cite{FV2} for the
function $u_m$. Finally, we state the 
self-adjointness of Macdonald-Ruijsenaars operators, and 
the qKZB heat equation with integration over a real cycle. 
The rest of the paper is devoted to 
the proof of these statements. Namely, in Section 3 we recall from 
\cite{EV3} the definition of the dynamical Weyl group and the 
dynamical Weyl group symmetry for trace functions; this fact
is vital for the proof of the main results. In Section 4
we prove the orthogonality relation and the qKZB-heat equation 
for traces in finite dimensional modules, using the techniques
of the paper \cite{EKi1}; besides of their independent interest, 
these results are used in the proof of the main results of Section 2. 
In Section 5, we show that the integrals considered in Section 2
are independent on the Weyl chamber in which the cycle of integration 
is situated; this is an important technical ingredient in the proof. 
In Section 6 we prove the orthogonality relation. In Section 7 we prove the 
qKZB heat equation. In Section 8 we prove the properties of
integral transforms. Finally, in Section 9 we prove the
self-adjointness of Macdonald-Ruijsenaars operators, and 
the qKZB heat equation with integration over a real cycle. 

We expect that the techniques of this paper can be extended 
to the case when a finite dimensional Lie algebra ${\frak g}$
is replaced with the affine Lie algebra $\widehat{\frak g}$, 
similarly to how the results of \cite{EV2} were extended to the affine case 
in \cite{ESV}. This would allow us to prove the orthogonality 
and qKZB-heat equation for trace functions of quantum affine algebras, and 
finally show (at least for $m=1$) that the function $u_m$ from \cite{FV2}
is  the trace 
function for $U_q(\widehat{sl_2})$ (up to renormalizations). This is a 
subject of future research. 

{\bf Acknowledgments.} It is our joy to dedicate this paper 
to the 90-th birthday of Izrail Moiseevich Gelfand.
He played a major role in the mathematical lives of both
authors. He also pioneered the idea to study special functions
by means of representation theory, of which this paper is an
example. 

The first author's work was partially 
supported by the NSF grant DMS-9988796, and done in part 
for the Clay Mathematics Institute. The second author's work was 
partially supported by the NSF grant DMS-9801582. 
The authors thank S. Tyurina for help in preparation of the paper.
\section{ Trace functions for $U_q(\g)$}

\subsection{Lie algebras and quantum groups}
Let $\g$ be a simple Lie algebra over $\C$ with root space decomposition 
$\g = \h \oplus ( \oplus_{\al \in \Si}\, \g_\al )$ where $\Si \subset \h^*$ is 
the set of roots.

Fix a system of simple roots $\al_1, ... , \al_r$. 
Let $\Si_\pm$ be  the set of positive
(negative) roots. Let $\n_{\pm}=\oplus_{\al\in \Si_{\pm}}\g_\al$. Then
$\g=\n_+\oplus\h\oplus\n_-$.

Let $(\,,\,)$ be
an invariant bilinear form on $\g$.
The form gives rise to a natural identification
$\h\to\h^*$, which we will sometimes use 
to make no distinction between $\h$ and $\h^*$.
This identification allows us to define a scalar product on
$\h^*$. We assume that the form is 
normalized so that $(\al, \al)=2$ for short roots.
We  use the same notation $(\,,\,)$ for the pairing $\h\T \h^*\to\C$.

We denote by 
$\h_\R$  (resp. $ \hr^*$ ) the real subspace of $\h$ (resp. $\h^*$).
The form $(\,,\,)$ is positive definite on $\hr$.

We use the  notation: 

$Q {} = \oplus_{i=1}^r \Z 
\alpha_i$ - the root lattice;

$Q_+=\oplus_{i=1}^r\Z_{\ge 0}\alpha_i$;

$Q^\vee=\oplus_{i=1}^r\Z\al_i^\vee$ - the dual root lattice,
where $\al^\vee=2\al/(\al,\al)$;

$P {} =\{\la\in\h\,|\, (\la,\al^\vee_i)\in\Z\}$ - the weight lattice;

$P_+=\{\la\in\h\,|\, (\la,\al^\vee_i)\in\Z_{\ge 0}\}$ - the 
cone of dominant integral weights;

$\om_i\in P_+$ - fundamental weights: $(\om_i,\al^\vee_j)=\dl_{ij}$;

$\rho {} ={1\over 2}\sum_{\al\in\Si_+}\al=\sum_{i=1}^r\om_i$;

Define a partial order on $\h$  putting $\mu<\la$ if $\la-\mu\in Q_+$.

A vector $
\la$ in $\hr$ or $\hr^*$ will be called {\it big} if
 $|\, (\la, \al_i)\,|\, >> \,0$ for $i=1,..., r$.
More precisely, a statement holds for big $\lambda$ 
if there exists a constant $K$ (possibly depending on some
previously fixed data) such that it holds for $\lambda$ 
satisfying the inequalities $|(\lambda,\alpha_i)|>K$ for all $i$. 

Let $s_i:\h^*\to\h^*$ denote 
the i-th simple reflection, defined by 
$$
s_i\ \la \ =\ \la-(\al_i^\vee,\la) \, \al_i\,.
$$ 
Let $\W$ be the Weyl group, generated by $s_1,...,s_r$.
%The following relations are defining:
%\bean\label{rela}
%s_i^2=1, \qquad (s_is_j)^m=1 \qquad \text{for}\qquad m=2,3,4,6,
%\notag
%\eean
%where $m=2$ if $\al_i$ and $\al_j$ are not neighboring in $\Gamma$,
%otherwise, $m=3,4,6$ if 1,2,3 lines respectively connect $\al_i$ and $\al_j$ in $\Gamma$.
For an element $w\in \W$,  denote by $l(w)$ 
the length of the minimal (reduced) presentation of $w$ as a product of generators
$s_1,...,s_r$.

We will also consider the ``dot'' action of the Weyl group on $\h^*$  defined by
$$
w \cdot v = w(v+\rho) - \rho.
$$
%where $\rho = \frac 12 \sum_{\alpha \in \Delta^+} \alpha.$

Let $(a_{ij})$ be the Cartan matrix of $\g$. 
Let $d_i$ be the relatively prime positive integers such that $(d_ia_{ij})$ 
is a symmetric matrix. Let 
$e_i,f_i,h_i$ be the Chevalley generators of $\g$.

Let $\kappa$ be a purely imaginary number, Im $\kappa \, < 0$, 
and $q=e^{{\pi i\over \kappa}}$, $0 < q < 1$. For any operator 
$A$, we denote $e^{{ \pi i\over \kappa} A}$ by $q^A$. 

Let $U_q(\g)$  be the Drinfeld-Jimbo quantum group corresponding 
to $\g$. Namely, $U_q(\g)$ is a Hopf algebra 
with generators
$E_i,F_i$, $i=1,\dots,r$, $q^h$, $h\in\h$, with relations:
$$
q^{x+y}=q^xq^y\ \text{for} \, x,y\in \h, q^0=1,\qquad 
q^hE_jq^{-h}=q^{\alpha_j(h)}E_i,\qquad q^hF_jq^{-h}=q^{-\alpha_j(h)}F_i\, ,
$$
$$E_iF_j-F_jE_i=\delta_{ij}\frac{q^{d_ih_i}-q^{-d_ih_i}}{q^{d_i}-q^{-d_i}},$$
$$\sum_{k=0}^{1-a_{ij}} (-1)^k \bmatrix 1-a_{ij} \\ k \endbmatrix_{q_i}
E_i^{1-a_{ij}-k}E_jE_i^k=0,\quad i\ne j,$$
$$\sum_{k=0}^{1-a_{ij}} (-1)^k \bmatrix 1-a_{ij} \\ k \endbmatrix_{q_i}
F_i^{1-a_{ij}-k}F_jF_i^k=0,\quad i\ne j.$$
where $q_i = q^{d_i}$ and we used the notation
$$ 
\bmatrix n \\ k \endbmatrix _q=\frac{[n]_q!}{[k]_q! [n-k]_q!},
\quad  [n]_q! = [1]_q \cdot [2]_q \cdot \dots \cdot [n]_q,
\quad  [n]_q= \frac {q^n - q^{-n}}{q-q^{-1}}\,.
$$

The comultiplication $\Delta,$ antipode $S$, and counit $\epsilon$
in $U_q(\g)$ are given by
$$
\Delta(E_i) = E_i\otimes q^{d_ih_i} + 1\otimes E_i, \quad
  \Delta(F_i) = F_i\otimes 1 + q^{-d_ih_i}\otimes F_i, \quad
  \Delta(q^h) = q^h \otimes q^h\,,
$$
$$
S(E_i)=-E_iq^{-d_ih_i},\quad S(F_i)=-q^{d_ih_i}F_i,\quad S(q^h)=q^{-h}\,,
$$
$$
\epsilon(E_i) = \epsilon(F_i) = 0,\quad \epsilon(q^h) = 1
\,.
$$

\subsection{ Intertwiners and trace functions}
Let $M_\mu$ be the Verma module over $U_q(\g)$ with highest weight 
$\mu$ and  highest weight vector  $x_\mu$. 
Let $V$ be a finite dimensional representation of $U_q(\g)$, 
and $v\in V$ a vector of weight $\mu_v$. It is well known 
that for generic $\mu$ there exists a unique intertwining 
operator $\Phi_\mu^v:M_\mu\to M_{\mu-\mu_v}\otimes V$ such that 
$\Phi_\mu^vx_\mu=x_{\mu-\mu_v}\otimes v+l.o.t.$ 
(here l.o.t. denotes the lower order 
terms, i.e. the terms of smaller weight in the first
component). 
It is useful to consider the ``generating function'' 
of such operators, $\Phi_\mu^V\in \Hom_\C(M_\mu,\oplus_\mu M_\nu\otimes V\otimes V^*)$, 
defined by 
$$
\Phi_\mu^V=\sum_{v\in B}\Phi_\mu^v\otimes v^*,
$$
 where the summation is 
over a homogeneous basis $B$ of $V$, and
$v^*$ are elements of the dual basis.

For $v\in V$ consider 
\bq\label{tr}
\Psi^{v}(\la,\mu)=
\text{Tr}|_{M_\mu}
(\Phi^{v}_{\mu} q^{2\la})\,, 
\end{equation}
a formal power series in
$V[0]\otimes q^{2(\la,\mu)}\C 
[[q^{-2(\la,\alpha_1)},...,q^{-2(\la,\alpha_r)}]]$.
This series converges (in a suitable region of values 
of  parameters) to a function
of the form $q^{2(\la,\mu)}f(\la,\mu)$, where $f$ is a rational function 
in $q^{2(\la,\alpha_i)}$ and $q^{2(\mu,\alpha_i)}$,
which is a finite sum of products of functions of $\la$ and
functions of $\mu$ (\cite{ESt},\cite{EV2}). This function %\Ref{tr} 
is called {\it the trace function}. 

The {\it universal trace function}
with values in $V[0]\otimes V^*[0]$ is the function
\bea\label{}
\Psi^{V}(\la,\mu)=\sum_{v\in B} 
\Psi^{v}(\la,\mu)\otimes v^*\,.
\eea
We have 
$\Psi^{V}(\la,\mu)=
\text{Tr}(\Phi^{V}_{\mu}q^{2\la})$.
We will consider the universal trace function as a function of $\la, \mu$ with values in 
$\End (V[0])$.

{\bf Example 1.}  
Let $\g=sl_2$. In this case we represent weights by complex
numbers, so that the unique fundamental weight corresponds to
$1$. Let $V=L_2$ be the irreducible 3-dimensional 
representation. Then 
$$
\Psi^V(\la,\mu)\ =\ \frac{q^{\la\mu}}{1-q^{-2\la}} \left(
1+(q^2-q^{-2})\frac{q^{-2\la}}{(1-q^{2\mu})(1-q^{-2(\la-1)})}\right)\,.
$$
(Since $V[0]$ is 1-dimensional, we view 
 $\Psi_V$ as a scalar function).

Let $V,W$ be finite dimensional representations of $U_q(\g)$. 
The fusion matrix is the operator $J_{WV}(\mu):W\otimes V\to W\otimes V$ defined 
by the formula
$$
(\Phi_{\mu-\mu_v}^{w}\otimes 1)\ \Phi_{\mu}^{v}\ =\
\Phi_\mu^{J_{WV}(\mu)(w\otimes v)},
$$
see \cite{EV1}. The exchange matrix 
$R_{VW}(\lambda)$ is defined by the formula 
$$
R_{VW}(\lambda):=J_{VW}(\lambda)^{-1}{\mathcal
R}^{21}_{WV}J^{21}_{WV}(\lambda),
$$
where ${\mathcal R}_{WV}$ is the R-matrix acting on $W\otimes V$. 

The universal 
fusion matrix $J(\la)$ takes values in a completion of 
$U_q(\g)\otimes U_q(\g)$ and
 gives $J_{VW}(\la)$ when evaluated in representations $V$, $W$.
%cf. also \cite{ABRR}, \cite{JKOS}. 
The universal exchange matrix is defined by the formula 
$R(\lambda):=J(\lambda)^{-1}{\mathcal R}^{21}J^{21}(\lambda)$,
where ${\mathcal R}$ is the universal R-matrix of $U_q(\g)$. 

If $J(\la)=\sum_i a_i\otimes b_i$ is the universal fusion matrix,
set $Q(\la)=\sum_i S^{-1}(b_i)a_i$, where $S$ is the antipode. This
sum defines an operator $Q_V(\la):V\to V$ invertible for generic
$\la$.

Let 
$$
\delta_q(\la)\ =\ \prod_{\alpha>0}\ (q^{(\la,\alpha)}-q^{-(\la,\alpha)})
$$
 be the Weyl denominator.

Introduce the renormalized trace function
 $$
 F^{V}(\la,\mu)\ =\ \delta_q(\la)\ \Psi^V(\la,-\mu-\rho)\ Q^{-1}_V(-\mu-\rho)\ .
$$

Let $V^*$ be the space dual to $V$ with the $U_q(\g)$ module structure defined by
the antipode. We have the following symmetry property \cite{EV2},
\bean\label{Sym}
F^{V^*}(\la,\mu)^*\ = \ F^V(\mu,\la)
\eean
where the 
values of both functions are regarded as linear operators on $V[0]$.

The following lemma describes the location of poles of
$F^V(\lambda,\mu)$. 

\begin{lemma}\label{simppol} The poles of the function $F^V(\lambda,\mu)$ with
respect to $\lambda$ are simple, and the divisor of poles 
is contained in the union of hyperplanes defined by
$[(\lambda,\alpha)-k(\alpha,\alpha)/2]_{q}=0$, where $k=1,2,
\ldots,N$, 
with $N=N(V),$
and $\alpha\in \Sigma_+$. 
\end{lemma}

\begin{proof} The lemma follows from Proposition 6.3 in \cite{EV2}.
\end{proof} 

 {\bf Example 2.} 
 Let $\g=sl_2$  and let $V=L_2$ be the 3-dimensional irreducible
 representation. Then
 $$
 F^V(\la,\mu)\ =\
 q^{-\la\mu}\ \frac{q^{2(\la+\mu)}-q^{2\la-2}-q^{2\mu-2}+1}{(1-q^{2\la-2})
 (1-q^{2\mu-2})} \ .
 $$
 Notice that the renormalized trace function is symmetric
 in $\la$ and $\mu$.

More generally, let $V=L_{2m}$.
Then according to \cite{EV2}, formula (7.20),
 \bea
F^V(\lambda,\mu) & = & q^{-\lambda\mu}\prod_{j=1}^{m}
\frac{q^{-2\mu-2j}-1}{q^{-2\mu-2j+2}-q^{-2m}}
\eea
\bea
&\times &
\sum_{l=0}^mq^{2m+l(l-1)/2}(q-q^{-1})^l
\frac{[m+l]_q!}{[l]_q![m-l]_q!}
\frac{q^{-2l\lambda}}{\prod_{j=1}^{l}(1-q^{-2(\mu+j)})
\prod_{j=1}^l(1-q^{-2(\lambda-j)})},
\eea
This function is symmetric in $\lambda$ and $\mu$, although
it is not obvious from the formula.

{\bf Example 3.} Let $\g=sl_n$, $V=S^{kn}\Bbb C^n$.
In this case, $V[0]$ is 1-dimensional, and $F^V(\lambda,\mu)$
is the Macdonald function studied in \cite{ESt} and \cite{Cha}. 

\subsection{Main results}
For $\xi \in \hr^*$ consider the imaginary subspace 
$C_\xi\, =\, \xi\,+\,i \hr^*$ in
$\h^*$. This subspace is invariant with respect to translations by
$\kappa Q^\vee$. The quotient $C_\xi/\kappa Q^\vee$ is a torus.
The theorems formulated below involve integration of functions over 
$C_\xi$ and over the torus $C_\xi/\kappa Q^\vee$. 
The integration will be performed with respect to translation invariant 
measures. Abusing notation, we will denote translation invariant
measures on $C_\xi$ and $C_\xi/\kappa Q^\vee$ by the same symbol 
$d\lambda$. The normalization of these measures is as follows. 
The measure $d\la$ on $C_\xi/\kappa Q^\vee$
is normalized by the condition 
$$
\int_{\ \  C_\xi / {\kappa Q^\vee}  }\  \, d\la\,=\, 1,
$$
while the measure $d\lambda$ on $C_\xi$ 
is normalized so that
$$
\ \int_{C_\xi}\ q^ {-(\la,\la)}\ d\la=1\ .
$$
This agreement is convenient because 
it allows one to get rid of normalization constants.
It will be kept throughout the paper. 
We warn the reader that with these normalizations, the 
direct image of the measure $d\lambda$ on $C_\xi$ is not equal
to the measure $d\lambda$ on $C_\xi/\kappa Q^\vee$, but 
is only proportional to it with some proportionality coefficient 
$C$, whose value is easy to compute but irrelevant to us.  

We now state the main results of this paper.  
Notice that the function 
$q^{2(\la, \mu)}\ F^V(\la, \mu)$, considered as a 
function of $\la$ 
( or as a function of $\mu$) and being restricted to $C_\xi$,  
is ${\kappa Q^\vee}$ periodic.

\begin{thm}\label{1} (Orthogonality)

Let $V$ be a finite dimensional $U_q(\g)$ module.
Assume that $\mu, \nu, \xi \in \hr^*$ are such that $\mu - \nu$ belongs to the weight lattice
$P$ and $\xi$ is big. 
Then
\bean
\int _{ C_\xi / {\kappa Q^\vee} }\, F^V(\mu, -\la)\ F^V(\la, \nu)\ d\la\ =\
\dl_{\mu, \nu}\ Q^{-1}_V (-\mu - \rho) \ .
\notag
\eean
\end{thm}

\begin{thm}\label{2} (qKZB-heat Equation)

Let $V$ be a finite dimensional $U_q(\g)$ module.
Assume that $\mu, \nu, \xi$ lie in $ \hr^*$,  and $\xi$ is big.
Then
\bea
\int _{ C_\xi  }\ F^V(\mu, -\la)\ F^V(\la, \nu)\ q^{-(\la, \la)}\ d\la\ 
=\ \ q^{(\mu,\mu) + (\nu, \nu)}\ F^V(\mu, \nu).\
\notag
\eea
\end{thm}

{\bf Remark 1.} If $V$ is the trivial representation, then 
$F^V(\lambda,\mu)=q^{-2(\lambda,\mu)}$, and 
Theorems \ref{1} and \ref{2} are obvious.

{\bf Remark 2.} 
 For $\g=sl_2$ the trace functions are given explicitly in
Example 2 of Section 2.2. 
Even in that case, Theorem \ref{1} and especially Theorem \ref{2}
are rather nontrivial integral identities.
More specifically, as seen from Example 2 above, 
the function $F^V$ for $V=L_{2m}$ is
a sum of $m+1$ products. So the left hand side of formulas in
Theorems \ref{1} and \ref{2} is a sum of $(m+1)^2$ integrals,
while the right hand side is a single product (or zero) in Theorem
\ref{1}, and a sum of $m+1$ products in Theorem \ref{2}.
Moreover, a careful computation shows that each of the 
$(m+1)^2$ individual integrals on the left hand side of Theorem \ref{2}
is non-elementary, and an elementary answer on the right hand side 
is obtained only as a result of cancellation.
\subsection{Integral transforms}

For $\xi,\eta \in \hr^*$ consider the imaginary subspace $C_\xi\, =\, \xi\,+\,i\hr^*$ in
$\h^*$ and the real subspace $D_\eta\, =\, i \eta\,+\,\hr^*$.

Let ${\mathcal S}(C_\xi)$ 
and ${\mathcal S}(D_\eta)$ be the Schwartz spaces of functions on
$C_\xi$ and 
$D_\eta$ 
respectively. 
Introduce the spaces
${\mathcal S}_\eta(C_\xi)=\{\phi: C_\xi\to \Bbb C\ |\
q^{-2i(\eta,\lambda)}
\phi(\lambda)
\in S(C_\xi)\},$
${\mathcal S}_{\xi}(D_\eta)=\{\phi: D_\eta\to \Bbb C 
|q^{2(\xi,\mu)}\phi(\mu)\in S(E_\eta)\}$.
Obviously, these spaces are canonically isomorphic 
to ${\mathcal S}_\eta(C_\xi),{\mathcal S}_\xi(D_\eta)$.
The modified Fourier transform $f(\lambda)\mapsto \hat
f(\mu):=\int_{C_\xi}q^{2(\lambda,\mu)}f(\lambda)d\lambda$
defines an isomorphism ${\mathcal S}_\eta(C_\xi)\to {\mathcal S}_\xi(D_\eta)$.
The inverse transform $g(\mu)\to g^\vee(\lambda)$ is given by
the formula $g^\vee(\lambda)=\int_{D_\eta}q^{-2(\lambda,\mu)}g(\mu)d\mu$. 
This fixes uniquely a normalization of the Lebesgue measure
$d\mu$ on $D_\eta$, which will be used from now on. 

For a finite dimensional $U_q(\g)$ module $V$ consider the function $F^V(\la, \mu)$.
Consider two integral transformations
\bean\label{Im}
K^V_{\IM}\ : \ {\mathcal S}_\eta(C_\xi)\otimes V[0] \ &\to & \ 
{\mathcal S}_\xi(D_\eta) \otimes V[0], 
\\
f(\la) \ & \mapsto & \int_{C_\xi}\ F^V(\mu, -\la)\ f(\la) \ d \la\ ,
\notag
\eean
and
\bean\label{Im1}
K^V_{\RE}\ : \ {\mathcal S}_\xi(D_\eta)\otimes V[0] \ &\to & \ {\mathcal S}_\eta(C_\xi)\otimes V[0], 
\\
f(\la) \ &\mapsto & \int_{D_\eta}\ F^V(\la, \mu)\ Q(-\mu-\rho) \ f(\mu) \ d \mu\ 
\notag
\eean

\begin{thm} \label{transform}

Assume that $\xi\in \h_{\Bbb R}$ is big and $\eta\in \h_{\Bbb R}$
is generic. Then the integral transforms are well defined, 
continuous in the Schwartz topology, and are inverse to each other,
$$
K^V_{\IM} \ K^V_{\RE} \ = \text{Id}\ ,
\qquad
K^V_{\RE} \ K^V_{\IM} \ = \text{Id}\ .
$$
\end{thm}

The proof of Theorem \ref{transform} occupies Section 8. 

\subsection{Self-adjointness of the Macdonald-Ruijsenaars
operators}

One of the basic facts of Macdonald's theory is that Macdonald's
operators are self-adjoint. It turns out that a similar statement holds
for the Macdonald-Ruijsenaars operators introduced in
\cite{EV2}. 

Namely, consider the  algebra ${\bold D}$ of scalar difference
operators. 
This algebra is generated by 
meromorphic functions $y(\lambda)$ on $\h^*$ and elements 
$T_\beta$ with defining relation
\bea
T_\beta y(\lambda)=y(\lambda+\beta)T_\beta.
\eea
The algebra $\bold D$ acts 
the space of functions on $\h^*$: functions act by
multiplication, and 
$T_\beta f(\lambda)=f(\lambda+\beta)$. 

Let $V$ be a finite dimensional $U_q(\g)$-module, and 
consider the algebra ${\bold D}_V={\bold D}\otimes {\rm End} V[0].$
Define an anti-homomorphism 
${\bold D}_V \to {\bold D}_{V^*}$ as follows:  
$\phi(\lambda)^*=\phi(-\lambda)$ for scalar functions,
$T^*_\beta=T_\beta$, and
$A\to A^*$ for $A \in {\rm End} V^*[0].$
It is obvious that for $L\in {\bold D}_V$, the operator $L^*$ is the formal
adjoint of $L$ with respect to the inner product $\langle f,g\rangle:=\int
(f(\lambda),g(-\lambda))d\lambda$. (We are not specifying the
contour of integration since it is not important for computing
the formal adjoint). 

Let ${\mathcal D}_{U,V}\in {\bold D}_V$ be 
the Macdonald-Ruijsenaars operator 
corresponding to the representation $U$ of $U_q(\g)$
(see \cite{EV2}). Namely, ${\mathcal D}_{U,V}$ is the 
difference operator acting on functions of $\lambda$ with values
in $V[0]$, which is defined by the formula
\bea
({\mathcal D}_{U,V}f)(\lambda)=\sum_\nu
\Tr|_{U[\nu]}(R_{UV}(-\lambda-\rho))
f(\lambda+\nu).
\eea 
According to \cite{EV2}, these operators commute for different
$U$ (i.e, form a quantum integrable system), and 
the trace function $F^V(\lambda,\mu)$ is their
common eigenfunction. 

\begin{thm}\label{selfad} One has
${\mathcal D}_{U,V}^*={\mathcal D}_{U,V^*}$
\end{thm}

\subsection{qKZB-heat equation with integration over a real cycle.}

\begin{thm}\label{realint} One has
\bea
\int_{D_\eta}F^V(\lambda,\mu)q^{(\mu,\mu)}Q(-\mu-\rho)
F^V(\mu,\nu)d\mu=q^{-(\lambda,\lambda)-(\nu,\nu)}
F^V(\lambda,\nu).
\eea
\end{thm}

The proof of this theorem is given in Section 9.

\section {Trace Functions and the Dynamical Weyl Group}
\subsection{The dynamical Weyl group, \cite {TV, EV3, STV}}
Recall that a nonzero vector in a $U_q(\g)$-module
is said to be singular if it is annihilated by $E_i$ for all $i$.

Let  $w=s_{i_1}\ldots s_{i_l}$ be a reduced
decomposition of $w\in \W$. Set $\al^{l}=\al_{i_l}$ and $\al^{j}=(s_{i_l}\ldots
s_{i_{j+1}}) (\al_{i_j})$ for $j=1,\ldots,l-1$. For $\mu \in \h^*$ let
$n_j=2\frac{(\mu+\rho,\alpha^j)} {(\alpha^j,\alpha^j)}$. For a
dominant $\mu \in P_+$, the numbers $n_j$ are positive integers. Let
$d^j=d_{i_j}$ (where $d_i$ are the symmetrizing numbers).
It is known that the collection of pairs of integers $(n_1,d^1),..., (n_k,d^k)$ and
the product $f_{\al_{i_1}}^{n_1}\cdots f_{\al_{i_l}}^{n_l}$ do
 not depend on the reduced decomposition.

 Define a vector
$x_{w\cdot\mu} \in M_\mu$ by 
\bean
x_{w\cdot\mu}\,=\, \frac{f_{\al_{i_1}}^{n_1}}{
[n_1]_{q^{d^1}}!} \ldots \frac{f_{\al_{i_l}}^{n_l}}{
[n_l]_{q^{d^l}}!}\,x_\mu\,.
\eean 
This vector is singular. It does not depend on the reduced decomposition.

Let $V$ be a finite dimensional $U_q(\g)$ module, and $w\in \W$.
According to \cite {TV, EV3, STV} there exists a unique operator $A_{w,V}(\mu)
\in \End (V)$ which rationally depends on 
$q^{2(\mu,\alpha_i)}$ and has the following properties.

Let $\mu \in P_+$ be a big vector.
Let $u\in V[\nu]$ for some $\nu\in \h^*$. Then
\begin{equation}
\Phi_\mu^u \ x_{w\cdot(\mu)}\  =
\ x_{w\cdot (\mu - \nu)} \ \otimes \ A_{w,V}(\mu)\ u \ +\ l.o.t.
\end{equation}
The collection of operators $\{ A_{w,V}(\mu)\}_{w\in\W}$ is 
called the dynamical Weyl group.

The operators of the dynamical Weyl group preserve the weight
decomposition of 
$V$ and satisfy
the cocycle condition. Namely, if  $w_1,w_2\in\W$, $l(w_1w_2)=l(w_1)+l(w_2)$, then
\bean\label{cocycle}
A_{w_1w_2,V}(\mu)\ =\ A_{w_1,V}(w_2\cdot\mu)\ A_{w_2,V}(\mu)\,.
\eean
Moreover, according to \cite{EV3}, 
on the subspace $V[0]$ this equation is satisfied without
the assumption $l(w_1w_2)=l(w_1)+l(w_2)$. 

Let ${\mathcal
A}_{w,V}(\lambda):=A_{w,V}(-\lambda-\rho)$. 
The trace function $F^V(\la, \mu)$ has the following symmetry property with respect to the 
dynamical Weyl group \cite{EV3},
\bean\label{symmetry}
F^V(\la, \mu) \ = 
\ {\mathcal A}_{w,V}(w^{-1}\la) \
F^V(w^{-1}\la, w^{-1} \mu)\  {\mathcal A}_{w,V^*}(w^{-1}\mu) ^*
\eean
for any $w\in \W$.

\subsection{Intertwiners of finite dimensional modules}
Let $L_\mu$ be
an irreducible finite dimensional representation with 
of $U_q(\g)$ 
with highest weight $\mu$. 
The intertwining operator $\Phi_\mu^v:  M_\mu\to M_\mu \otimes V$,
descends to an operator $\bar\Phi_\mu^v: L_\mu\to L_\mu \otimes V$.
We define the corresponding trace functions by
$$
\Psi^v_\mu(\la)\ =
\  \text{Tr}|_{L_\mu}(\bar\Phi_\mu^{v}
q^{2\la}), \qquad
\Psi^V_\mu(\la)\ =
\ \sum_{v\in B} \Psi^{v}_\mu(\la)\otimes v^*
$$
and the corresponding renormalized trace function as
$$
F^V_\mu(\lambda) \ =\ \delta_q(\lambda)\ \Psi^V_{-\mu-\rho}
(\lambda)\ Q_V^{-1}(-\mu -\rho)\,.
$$

If we regard universal trace functions  as linear operators on $V[0]$,
then we have the following generalized Weyl character formulas \cite{EV3}.
For a big dominant integral weight $\mu$, we have
$$
\Psi^V_\mu(\la) \ = \ \sum_{w \in \W} \ (-1)^w \ 
\Psi^V(\la, w\cdot \mu)\  A_{w,V}(\mu),
$$
and for a big anti-dominant integral weight $\mu$ we have
\bean\label{dWeyl}
F^V_\mu(\la) \ = \ \sum_{w \in \W} \ (-1)^w \ F^V(\la, w \mu) 
 \ ({\mathcal A}_{w, V^*}(\mu)^{-1})^*.
\eean
Here $(-1)^w$ denotes the sign of the element $w$.
The terminology is motivated by the fact 
that for $V=\Bbb C$ these formulas reduce to the usual Weyl
character formula. 

\section{Orthogonality and the 
qKZB-heat equation for finite dimensional modules}
\subsection{Statement of results}

\begin{thm}\label{1'} (Orthogonality)

Let $\mu, \nu$ be big dominant integral weights.
% such that 
%$|\, (\mu, \al_i)\,|\, >> \,0$ and $|\, (\nu, \al_i)\,|\, >> \,0$ for $i=1,..., r$. 
Let $V$ be a finite dimensional $U_q(\g)$ module, $V^*$ the dual module.
Let $v\in V[0]$ and $v_*\in V^*[0]$ be arbitrary vectors.
Consider the intertwiners $\bar\Phi_\mu^v: L_\mu\to L_\mu \otimes V$,
$\bar\Phi_\nu ^{v_*}: L_\nu \to L_\nu \otimes V^*$ and the corresponding trace 
functions $\Psi^v_\mu(\la), \ \Psi^{v_*}_\nu(\la)$. Let $\xi\in\hr^*$ be a big vector.
%such that
%$|\, (\xi, \al_i)\,|\, >> \,0$ for $i=1,..., r$. 
Then
\bea
{1\over | \W |}\ 
\int _{ C_\xi / {\kappa Q^\vee} }\ \delta_q(\la)\ \delta_q(-\la)
\ (\Psi^v_\mu(\la),  \Psi^{v_*}_\nu(- \la))\ d\la
\ = \ \delta_{\mu,\nu}\ ( Q_V(\mu)v, v_*)\,.
\notag
\eea
Here $|\W|$ is the number of elements in the Weyl group and $(\,.\,,\,.\,)$ is the 
pairing of vectors and covectors.
\end{thm}

\begin{thm}\label{2'} (qKZB-heat equation)

Under the assumptions of Theorem \ref{1'}, one has 
\bea
 {1\over | \W |}\ 
\int _{ C_\xi  }\ \delta_q(\la)\ \delta_q(-\la)\ 
 ( \Psi^v_\mu(\la),  \Psi^{v_*}_\nu(- \la))\ q^{-(\la, \la)}\ d \la 
\\ 
 = \ \delta_q(-\mu - \rho)\ q^{(\mu + \rho, \mu + \rho) + (\nu + \rho, \nu + \rho)}\
(Q_V(\mu) v, \Psi^{v_*}_\nu( - \mu - \rho))\,.
\notag
\eea
\end{thm}

\subsection{ Proof of Theorem \ref{1'}}

The function $\Psi_\nu^{v_*}(-\lambda)$ 
is the trace in $L_\nu$ of the operator $A:=\bar\Phi_\nu^{v_*}
q^{-2\lambda}$, acting from $L_\nu$ to $L_\nu\otimes V^*$. 
Therefore, the same function can be computed as the trace of the
dual operator $A^*$. The operator $A^*$ can be written as 
$q^{2\lambda}(\bar\Phi_\nu^{v_*})^*$. 

To interpret $(\bar\Phi_\nu^{v_*})^*$ as an
intertwiner, let us use the operation of left dual, $W\to {}^*W$,
on representations of $U_q(\g)$. Namely, ${}^*W$
is the usual dual of $W$ as a vector space, with 
the action of $U_q(\g)$ defined by the formula 
$\pi_{{}^*W}(a)=\pi_W(S^{-1}(a))^*$ (as opposed 
to $\pi_W(S(a))^*$ for the right dual $W^*$). 
Also, the intertwining operator 
$L_\nu\to Y\otimes L_\nu$ which sends the 
{\it lowest} weight vector $x^\vee_\nu$ to $y\otimes x^\vee_\nu+...$
for some $y\in Y$ will be denoted by $\tilde\Phi_\nu^y$.

With these definition, we see that $(\bar\Phi_\nu^{v_*})^*$ 
is an intertwiner $V\otimes {}^*L_\nu\to {}^*L_\nu$, 
which can also be viewed (by swapping $V$) as an intertwiner 
${}^*L_\nu\to {}^*V\otimes {}^*L_\nu$. 
Thus, $(\bar\Phi_\nu^{v_*})^*=\tilde\Phi_{\nu^*}^{v_*}$, where 
$\nu^*$ is the highest weight of ${}^*L_\nu$.
Thus, using the cyclic property of the trace, 
we find that $\Psi_\nu^{v_*}(-\lambda)=\Tr(\tilde\Phi_{\nu^*}^{v_*}q^{2\lambda})$. 

This shows that the expression
$(\Psi_\mu^v(\lambda),\Psi_\nu^{v_*}(-\lambda))$ 
can be represented in the form $\Tr|_{L_\mu\otimes
L_\nu^*}(X(q^{2\lambda}\otimes q^{2\lambda}))$, 
where $X\in {\rm End}(L_\mu\otimes
L_\nu^*)$ is the composition of the tensor product 
$\bar\Phi_\mu^v\otimes \tilde\Phi_\nu^{v_*}$ with the 
contraction $V\otimes {}^*V\to \Bbb C$. 

Now observe that $X$ is an intertwining operator.
Therefore, writing $L_\mu\otimes L_\nu^*$ as a direct sum 
$\oplus_\beta H_{\mu\nu^*}^\beta\otimes L_\beta$
(where $H_{\mu\nu^*}^\beta$ are the multiplicity spaces),
we can represent $X$ in the form 
$X=\oplus_\beta X_\beta\otimes 1_{L_\beta}$, 
where $X_\beta\in {\rm End}(H_{\mu\nu^*}^\beta)$. 
Hence, 
$$
(\Psi_\mu^v(\lambda),\Psi_\nu^{v_*}(-\lambda))=
\sum_\beta \Tr(X_\beta)\chi_\beta(q^{2\lambda}),
$$ 
where $\chi_\beta$ is the character of the representation
$L_\beta$. Therefore, by the Weyl orthogonality formula for
characters we get 
\bea
{1\over | \W |}\ 
\int _{ C_\xi / {\kappa Q^\vee} }\ \delta_q(\la)\ \delta_q(-\la)
\ (\Psi^v_\mu(\la),  \Psi^{v_*}_\nu(- \la))\ d\la
= \delta_{\mu,\nu}\Tr(X_0).
\notag
\eea
This immediately implies the theorem for $\mu\ne \nu$ 
(as in this case $H_{\mu\nu^*}^0=0$ and hence $\Tr(X_0)=0$). 
Thus, it remains to settle the case $\mu=\nu$. 

If $\mu=\nu$ then $H_{\mu\nu^*}^0$ is one dimensional, and $X_0$
is a number. So we need to compute this number. 

To compute $X_0$, pick $y\in L_\nu$, 
$f\in {}^*L_\nu$, and look at ${\rm cont}(X(y\otimes f))$, where 
we denote by ${\rm cont}: L_\nu\otimes {}^*L_\nu\to \Bbb C$ the contraction
operator. It is easy to see that \linebreak 
${\rm cont}(X(y\otimes f))=X_0f(y)$. On the other hand, 
it is easy to check from the definition of $X$ that 
${\rm cont}(X(y\otimes f))=af(y)$, where 
$a$ is found from the equation 
$$
(1\otimes (,))(\bar\Phi^{v_*}_\nu\otimes 1)\bar\Phi^v_\nu=a\cdot 1_{L_\nu}
$$
Thus, $X_0=a$, and it can be found in terms of the fusion
matrix. Namely,  we find
$$
(1\otimes (,))\bar\Phi^{J_{V^*,V}(\nu)(v_*\otimes v)}_\nu=X_0\cdot 1_{L_\nu},
$$
hence $X_0$ is the contraction of $J_{V^*,V}(\nu)(v_*\otimes v)$,
which equals $(S(Q)(\nu)v,v_*)$. However, as follows from
Proposition 2.13 of \cite{EV2}, $Q=S(Q)$ on the zero weight
subspace. This implies Theorem \ref{1'}. 

{\bf Remark.} Let $\g=sl_n$, and $V$ be the q-deformation of the
representation $S^{kn}\Bbb C^n$. 
In this case, $V[0]$ is one-dimensional (so $\Psi_\nu^v$ can be viewed
as a scalar function and is independent of $v$ up to scaling), and 
the weights $\nu$ for which the operators $\bar\Phi_\nu^v$ exist
are those of the form $\mu+k\rho$, where $\mu$ is a dominant
integral weight. Moreover, as was shown in \cite{EKi2}, 
the function $\Psi_{k\rho}(\lambda)$ is given by an explicit 
product formula, while $\Psi_{\mu+k\rho}(\lambda)/\Psi_{k\rho}(\lambda)$ is 
the Macdonald polynomial with highest weight $\mu$. Thus, 
using Proposition 41 from \cite{EV3} (the determinant formula for
$Q$), we obtain the Macdonald inner product identities for Macdonald's polynomials
of type $A$ (see e.g. \cite{EKi3} for the formulation). 
Another (more complicated) representation
theoretic proof of these identities was given in \cite{EKi3}.

\subsection { A remark on theta functions and the Kostant identity} 
Let $f$ be a smooth function on 
 $C_\xi$, which is periodic with respect to the lattice 
${\kappa Q^\vee}$. The function  $f$ can be decomposed into a Fourier series 
with respect to the basis $q^{2(\la,\beta)}$, $\beta\in P$. 

Let $\gamma(\la)$ denote the theta-function 
$$
\gamma(\la)\ =\ \sum_{\beta\in P}\ q^{\beta^2} \ q^{2(\la,\beta)} \ .
$$
The following lemma is standard. 

\begin{lemma} \label{theta} We have
$$
\int_{C_\xi}f(\la)\ q^{-\la^2}\ d\la\  =\ 
 \int_{C_\xi/{\kappa Q^\vee}}\ \ f(\la)\gamma(\la)\ d\la\,.
$$
\end{lemma} 

\begin{proof} We have
$$
\int_{C_\xi}\ f(\la)\ q^{-\la^2}\ d\la\ =\
C\int_{C_\xi/{\kappa Q^\vee}}\ f(\la)\ \sum_{\chi\in {\kappa Q^\vee}} \ q^{-(\la+\chi)^2} d\la\ .
$$
It follows form comparison of Fourier coefficients  that
$$
C\sum_{\chi\in {\kappa Q^\vee}}\  q^{-(\la+\chi)^2} \ =\ \gamma (\la)\ .
$$
Thus the lemma is proved.
\end{proof}

Recall that the quantum dimension, $\dim_qL_\nu$, of the
representation $L_\nu$, is the number $\chi_\nu(q^{2\rho})$. 
One has $\dim_q L_\nu=\prod_{\alpha\in
\Sigma_+}[(\alpha,\nu+\rho)]_q=\delta_q(\nu+\rho)/\delta_q(\rho)$. 

\begin{thm} (Kostant, \cite{Kos}) 
One has 
$$
\gamma(\lambda)=K\sum_{\beta\in
P_+}q^{(\beta,\beta+2\rho)}
\chi_\beta(q^{2\lambda})\dim_q(L_\beta), 
$$
where $K=\prod_{\alpha\in
\Sigma_+}(1-q^{2(\alpha,\rho)})$.
\end{thm} 

\subsection{ Proof of Theorem \ref{2'}}

The proof follows the ideas of \cite{EKi1}. 

By Lemma \ref{theta}, the statement of Theorem \ref{2'}
is equivalent to the equality
\bea
 {1\over | \W |}\ 
\int _{ C_\xi/{\kappa Q^\vee}  }\ \delta_q(\la)\ \delta_q(-\la)\ 
 ( \Psi^v_\mu(\la),  \Psi^{v_*}_\nu(- \la))\ \gamma(\la)\ d \la 
\\ 
 = \ \delta_q(-\mu - \rho)\ q^{(\mu + \rho, \mu + \rho) + (\nu + \rho, \nu + \rho)}\
(Q_V(\mu) v, \Psi^{v_*}_\nu( - \mu - \rho))\,.
\notag
\eea

Using the argument and notation of the proof of Theorem \ref{1'},
and also the Kostant identity, we can
rewrite the left hand side of this equation in the form
\bea
LHS= {K\over | \W |}\ 
\int _{ C_\xi/{\kappa Q^\vee}  }\ \delta_q(\la)\ \delta_q(-\la)\ 
\sum_{\beta\in P_+}\Tr(X_\beta)\chi_\beta(q^{2\lambda})\times
\\
\sum_{\beta'\in
P_+}q^{(\beta',\beta'+2\rho)}
\chi_{\beta'}(q^{2\lambda})\dim_q(L_{\beta'})d\la 
\notag
\eea
Thus by the Weyl orthogonality theorem for characters, 
we have 
\bea
LHS=  
K\sum_{\beta\in P_+}\Tr(X_\beta)q^{(\beta,\beta+2\rho)}
\dim_q(L_{\beta}). 
\notag
\eea
Indeed, the only nonzero contributions come from the case
$\beta'=\beta^*$, while the expressions $(\beta,\beta+2\rho)$ 
and $\dim_q(L_\beta)$ are invariant under the transformation
$\beta\to\beta^*$. 

Let $u$ be the Drinfeld element of $U_q(\g)$ \cite{Dr}. 
Namely, $u$ is an element of a completion of $U_q(\g)$ defined by
the formula $u=\sum S(b_i)a_i$, where $\sum a_i\otimes b_i$ is
the universal R-matrix $\mathcal R$ of $U_q(\g)$. Drinfeld showed that 
$u$ acts on $L_\beta$ as $q^{-(\beta,\beta+2\rho)}q^{2\rho}$, and that 
$uS(u)^{-1}=q^{4\rho}$. Thus, $S(u)^{-1}$ acts in $L_\beta$ by
$q^{(\beta,\beta+2\rho)}
q^{2\rho}$, and hence
$q^{(\beta,\beta+2\rho)}\dim_qL_\beta=\Tr|_{L_\beta}(S(u)^{-1})$. 
Therefore, we get 
\bea
LHS=  
K\sum_{\beta\in P_+}\Tr(X_\beta\otimes S(u)^{-1}|_{L_\beta})=
K\Tr(X\Delta(S(u)^{-1})). 
\notag
\eea
Now, as was also shown by Drinfeld \cite{Dr}, one has 
$\Delta(u)=({\mathcal R}^{21}{\mathcal R})^{-1}(u\otimes u)$. 
Hence, 
$\Delta(S(u)^{-1})={\mathcal R}^{21}{\mathcal R}(S(u)^{-1}\otimes
S(u)^{-1})$ (we use that $(S\otimes S)({\mathcal R})={\mathcal
R}$). 
Substituting this into the formula for the LHS, 
and remembering the definition of $X$, we get 
$$
LHS=K(,)\circ \Tr|_{L_\mu\otimes {}^*L_\nu}
\left((\bar\Phi_\mu^v\otimes \tilde\Phi_\nu^{v_*}){\mathcal
R}^{21}{\mathcal R}
(S(u)^{-1}\otimes S(u)^{-1})\right)
$$
(the trace takes values in $V\otimes {}^*V$, and $(,)$ denotes
the contraction $V\otimes {}^*V\to \Bbb C$). 

To compute this trace, let us look at the 
trace in one of the factors, i.e. 
\bea\label{tracesingle}
\Tr|_{{}^*L_\nu}((1\otimes \tilde\Phi_{\nu^*}^{v_*}) {\mathcal
R}^{21}{\mathcal R}
(1\otimes S(u)^{-1}))=
q^{(\nu,\nu+2\rho)}T,\\
T:=\Tr|_{{}^*L_\nu}((1\otimes \tilde\Phi_{\nu^*}^{v_*}) 
{\mathcal R}^{21}{\mathcal R} (1\otimes q^{2\rho})).
\eea
It is easy to show that if $Y,Z$ are $U_q(\g)$ modules, $\dim
Z<\infty$, and $\Phi\in {\rm End}_{U_q(\g)}(Y\otimes Z)$, then 
$\Tr|_Z(\Phi(1\otimes q^{2\rho}))$ is an intertwining operator
$Y\to Y$. Therefore, $T$ is an intertwiner $L_\mu\to L_\mu\otimes
{}^*V$, and the operator 
$$
B:=(1\otimes (,))
(\bar\Phi_{\mu}^{v}\otimes 1)T: L_\mu\to L_\mu
$$
is a scalar. 

Now, we have 
$$
K^{-1}q^{-(\mu,\mu+2\rho)-(\nu,\nu+2\rho)}LHS=\Tr(Bq^{2\rho})=
B\dim_q(L_\mu)= 
B\frac{\delta_q(\mu+\rho)}{\delta_q(\rho)}=
B\frac{\delta_q(-\mu-\rho)}{\delta_q(-\rho)},
$$
and 
$$
K=q^{2\rho^2}\delta_q(-\rho). 
$$
Therefore, after a simple calculation we find that 
the statement of Theorem \ref{2'}
is equivalent to the identity
\bean\label{B}
B=(Q_V(\mu)v,\Psi_\nu^{v_*}(-\mu-\rho)).
\eean
So it remains to prove formula (\ref{B}). 

To prove (\ref{B}), let us apply the operator $T$ to the highest
weight vector $x_\mu$ of $L_\mu$. We have ${\mathcal
R}(x_\mu\otimes y)=x_\mu\otimes q^\mu y$. Therefore, 
\bea
Tx_\mu=\Tr|_{{}^*L_\nu}((1\otimes \tilde\Phi_{\nu^*}^{v_*})
{\mathcal R}^{21} (1\otimes q^{\mu+2\rho}))x_\mu=\\
x_\mu\otimes \Tr|_{{}^*L_\nu}(\tilde\Phi_{\nu^*}^{v_*}q^{2(\mu+\rho)})+l.o.t
\notag
\eea
It follows from the proof
of Theorem \ref{1'} that 
$$
\Tr|_{{}^*L_\nu}(\tilde\Phi_{\nu^*}^{v_*}q^{2(\mu+\rho)})=
\Psi_\nu^{v_*}(-\mu-\rho).
$$
Therefore, 
$$
Tx_\mu=x_\mu\otimes \Psi_\nu^{v_*}(-\mu-\rho)+l.o.t.,
$$
i.e., $T=\bar\Phi_\mu^{\Psi_\nu^{v_*}(-\mu-\rho)}$. 
Thus, it follows form the definition of $B$ (as in the proof 
of Theorem \ref{1'}) that 
$$
B=(Q_V(\mu)v,\Psi_\nu^{v_*}(-\mu-\rho)).
$$
Theorem \ref{2'} is proved. 

{\bf Remark.} As we mentioned, in the case $\g=sl_n$,
$V=S^{kn}\Bbb C^n$, the functions $\Psi_\mu^v$ are, up to
normalization, Macdonald polynomials (of type A).   
Theorem \ref{2'} in this case coincides with Cherednik's
q-deformation of the Macdonald-Mehta identity \cite{Ch1}, 
and the proof we gave is the same as given in \cite{EKi1}. 

\subsection{Reformulation of Theorems \ref{1'} and \ref{2'}}

Let us now reformulate Theorems \ref{1'} and \ref{2'} 
in terms of renormalized trace functions $F_\mu^V(\lambda)$. 

\begin{thm}\label{1''} (Orthogonality)

Let $\mu, \nu$ be big dominant integral weights.
% such that 
%$|\, (\mu, \al_i)\,|\, >> \,0$ and $|\, (\nu, \al_i)\,|\, >> \,0$ for $i=1,..., r$. 
Let $V$ be a finite dimensional $U_q(\g)$ module, $V^*$ the dual module.
Consider the trace functions $F^V_\mu(\la), \ F^{V^*}_\nu(\la)$ as linear operators on
the corresponding zero weight subspaces.
 Let $\xi \in \hr^*$ be a big vector.
%such that
%$|\, (\xi, \al_i)\,|\, >> \,0$ for $i=1,..., r$. 
Then
\bea
{1\over | \W |}\ 
\int _{ C_\xi / {\kappa Q^\vee} }\ 
 F^{V^*}_\nu(- \la)^*  \ F^V_\mu( \la) \ d\la
\ = \ \delta_{\mu,\nu}\ Q_V^{-1}(-\mu-\rho)\,.
\notag
\eea
Here $ F^{V^*}_\nu(\la)^*$ denotes the operator on $V$ dual to the operator
$ F^{V^*}_\nu(\la)$ on $V^*$. 
\end{thm}

\begin{thm}\label{2''} (qKZB-heat equation)

Under the assumptions of Theorem \ref{1''}, one has
\bea
 {1\over | \W |}\ 
\int _{ C_\xi  }\ 
 F^{V^*}_\nu(- \la)^* \ F^V_\mu(\la) \ q^{- (\la, \la)}\ d\la 
\ = \ q^{(\mu, \mu) + (\nu, \nu)}\
 F^{V^*}_\nu(\mu)^* \,.
%\notag
\eea
\end{thm}

{\bf Corollary of Theorem \ref{2''}.} {\it Under 
the conditions of Theorem \ref{2''} we have}
\bean
 F^{V^*}_\nu(\mu)^* \ = \  F^{V}_\mu(\nu) \, .
\notag
\eean
Indeed, 
let us interchange $\mu$ and $\nu$ in the 
formula of Theorem \ref{2''},  change $\la$ to
$-\la$, and take the dual operators to the 
operators on the left and right hand sides 
of the formula. Then the left hand side will remain the same, while in the right hand side 
the operator $ F^{V^*}_\nu(\mu)^*$ will be replaced by $ F^{V}_\mu(\nu)$.

{\bf Remark.} For $\g=sl_n$, $V=S^{kn}\Bbb C^n$, 
this corollary reduces to the Macdonald's symmetry identity. 

\subsection{Proof of Theorem \ref{1''}}
 From the definition of $F^V_\mu(\lambda)$ we have
%$$
%F^V_\mu(\lambda) \ =\  \delta_q(\lambda) \ \Psi^V_{-\mu-\rho}(\lambda)\
%Q_V^{-1}(-\mu-\rho)\ .
%$$
%Hence
\bean\label{f}
\delta_q(\lambda) \ \Psi^v_{\mu}(\lambda)\ =\ F^V_{-\mu-\rho}
(\lambda)\  
Q_V (\mu) v\ .
\eean
Substituting this into Theorem \ref{1'} we get
$$
\frac{1}{|\W|} \ \int_{C_{\xi}/{\kappa Q^\vee}}\  (\ F^{V}_{-\mu-\rho}(\lambda)\ Q_V(\mu)v ,\
F^{V^*}_{-\nu-\rho}(-\lambda)\ Q_{V^*}(\nu)v_*\ )\  d\lambda \ =
\ \delta_{\mu \nu} \ (Q_V (\mu)v, v_*)\ .
$$
The integrand can be written as
$$
(Q_{V^*}(\nu)^* \ F^{V^*}_{-\nu-\rho}(-\lambda)^*\ F^{V}_{-\mu-\rho}(\lambda)\ Q_V(\mu)v,\  v_*) .
$$
By Proposition 2.13 in \cite{EV2}  we have $Q_{V^*}(\nu)^*|_{V[0]}=Q_{V}(\nu)|_{V[0]}$. Thus
$$
\frac{1}{|\W|} \ \int_{C_{\xi}/{\kappa Q^\vee}} 
\ (\ Q_{V}(\nu)\ F^{V^*}_{-\nu-\rho}(-\lambda)^*\ F^{V}_{-\mu-\rho}(\lambda)\ Q_V(\mu)v, \ v_*\ )
\ d\lambda \ =\
\delta_{\mu \nu} \ (Q_V (\mu)v, v_*)\ .
$$
That gives Theorem \ref{1''}.

\subsection{Proof of Theorem \ref{2''}}
Using \Ref{f} the formula of Theorem \ref{2'} can be written as
\bea
\frac{1}{|\W|} \ \int_{C_{\xi}} 
\ (\ Q_{V^*}(\nu)^*\ F^{V^*}_{-\nu-\rho}(-\lambda)^*\ F^{V}_{-\mu-\rho}(\lambda)\ Q_V(\mu)v, \ v_*\ )\
q^ {-(\la,\la)}\ d\la\ =
\\
\ q^{(\mu+\rho, \mu+\rho) + (\nu+\rho, \nu+\rho)}\ 
 \ (Q_{V^*}(\nu)^* \ F_{-\nu-\rho}^{V^*}(-\mu-\rho)^*\ Q_V (\mu)v, \  v_*)\ .
\eea
That gives Theorem \ref{2''}.

\section{Independence of  integrals on the choice of a Weyl chamber}

\subsection {Statement of the result}
For $\la\in \h^*$ write $\la = x + i y$ with $x, y \in\hr^*$.
Let $g(\la)$ be a holomorphic function on $\h^*$. We say that
$g(\la)$ is  rapidly 
decaying in the imaginary direction 
if for any positive $s,r$ there exists $C_{s,r}>0$
such that 
$$
|g(x+iy)|\le C_{s,r}(1+|y|)^{-s}
$$
as long as $|x|\le r$. 
An example of such a function is $g(\lambda)=q^{-(\lambda,\lambda)}$. 

Let $g(\la)$ be a holomorphic function on $\h^*$ 
which is invariant with respect to the 
standard Weyl group action and rapidly decaying in the imaginary direction. 
Let $\xi \in \hr^*$ be a big vector. Consider the integral
\bean\label{integral}
I(\xi)\ = \ \int_{C_\xi}\  g(\la)\ F^{V^*}(-\la, \mu)^*\ F^{V}(\la,\nu) \ d\la\ .
\eean
It is clear that the integral does not depend on the choice of
the big vector $\xi$ as long as the vector
$\xi$ belongs to the same Weyl chamber.

\begin{thm}\label{chamberindep}
Under the above conditions the integral $I( \xi )$
does not depend on the choice of the Weyl chamber containing $\xi$.
\end{thm} 

Theorem \ref{chamberindep} is proved in the next three subsections.

We will also need the following modification 
of Theorem \ref{chamberindep}.

Let $g(\la)$ be a holomorphic function on $\h^*/\kappa Q^\vee$ 
which is invariant with respect to the 
standard Weyl group action. 
Let $\xi \in \hr^*$ be a big vector. Consider the integral
\bean\label{integral1}
J(\xi)\ = \ \int_{C_\xi/\kappa Q^\vee}\  
g(\la)\ F^{V^*}(-\la, \mu)^*\ F^{V}(\la,\nu) \ d\la\ .
\eean
As before, it is clear that the integral does not depend on the choice of
the big vector $\xi$ as long as the vector
$\xi$ belongs to the same Weyl chamber.

\begin{thm}\label{chamberindep1}
Under the above conditions the integral $J( \xi )$
does not depend on the choice of the Weyl chamber containing $\xi$.
\end{thm} 

The proof of Theorem \ref{chamberindep1} is completely parallel
to the proof of Theorem \ref{chamberindep} and is omitted. 

\subsection{The $sl_2$ case}
First of all, let us prove Theorem \ref{chamberindep} for $\g=sl_2$ 
and $V=L_{2m}$, the irreducible representation 
with highest weight
$2m$ ($m\in \Bbb Z_+$). 
In this case, the theorem reduces to the equality 
$$
\int_{{\rm Re}\lambda=-a}
g(\lambda)F(-\lambda,\mu)F(\lambda,\nu)d\lambda=
\int_{{\rm Re}\lambda=a}
g(\lambda)F(-\lambda,\mu)F(\lambda,\nu)d\lambda,
$$
for large enough $a>0$, 
which (due to rapid decay of $g$ in the imaginary direction)
is equivalent to the statement that the sum of residues
of the integrand in the strip $|{\rm Re}(\lambda)|<a$ is zero. 

{\bf Remark.} Here for brevity we write $F$ instead of $F^V$. We
also do not distinguish between $F$ and $F^*$, since 
$V[0]$ is 1-dimensional. 
By Lemma \ref{simppol}, the poles of $F(\lambda,\mu)$ 
with respect to $\lambda$ are simple, and located at the points $1,2,...,m$
and their translates by $\kappa \Bbb Z$. 
Therefore, the poles of $G(\lambda):=F(-\lambda,\mu)F(\lambda,\nu)$
are simple as well and located at $\pm k+\kappa\Bbb Z$,
$k=1,...,m$. Since the function $G$ is quasi-periodic with period
$\kappa$ (i.e. it is multiplied by a constant under the shift by
$\kappa$), 
the cancellation of residues (and hence the theorem) follows from 

\begin{proposition}\label{resi} The residue of the function $G(\lambda)$ 
at $\lambda=k$ equals minus its residue at $\lambda=-k$. 
\end{proposition} 

The rest of the subsection is the proof of this proposition. 

Consider 
the function ${\tilde{\Psi}}(\lambda,\mu):=\Psi(\lambda,\mu)\prod_{i=0}^{m-1}
(q^{m-i}-q^{-m+i}).$

Using formula (7.19) for $Q(\mu)$ in \cite{EV2}, we have
\bea
F(\lambda,\mu)=(q^\lambda-q^{-\lambda})
\Psi(\lambda,-\mu-1)Q^{-1}(-\mu-1)=
\eea
\bea
(q^\lambda-q^{-\lambda})
{\tilde{\Psi}}(\lambda,-\mu-1)Q^{-1}(-\mu-1)
\prod_{i=0}^{m-1}\frac{1}{q^{-(\mu+1+i)}-q^{\mu+1+i}}=
\eea
\bea
(q^\lambda-q^{-\lambda})
{\tilde{\Psi}}(\lambda,-\mu-1)q^{2m}
\prod_{j=1}^{m}\frac{q^{-2\mu-2j}-1}{q^{-2\mu-2j+2}-q^{-2m}} 
\prod_{i=0}^{m-1}\frac{1}{q^{-(\mu+1+i)}-q^{\mu+1+i}}=
\eea
\bea
C(q)(q^\lambda - q^{-\lambda}){\tilde{\Psi}}(\lambda,-\mu-1)
\prod_{j=1}^{m}\frac{1}{q^{-\mu-j+1+m}-q^{\mu+j-1-m}}.
\eea
Changing $\mu$ to $-\mu$ in this formula, we get
\bea
F(\lambda,-\mu)=
C(q)(q^\lambda - q^{-\lambda}){\tilde{\Psi}}(\lambda,\mu-1) 
\prod_{j=1}^{m}\frac{1}{q^{\mu-j+1+m}-q^{-\mu+j-1-m}}
\eea

Now recall the following result. 

\begin{thm}\cite{ESt} \label{ESt} The function 
${\tilde{\Psi}}$ is holomorphic in $\mu$, 
and satisfies the quasi-invariance (or resonance) conditions
$$
{\tilde{\Psi}}(\lambda,k)={\tilde{\Psi}}(\lambda,-k-2),
\qquad
 k=1,...,m.
$$
\end{thm}

This result and the above formulas, after a short calculation,
lead to the identity
\bea
Res_{\mu=k}F(\lambda,\mu)=
\eea
\bea
F(\lambda,-k)\ln(q^{-2})\prod_{j=1,j\neq m+1-k}^{m}
\frac{1}{q^{-k-j+1+m}-q^{k+j-1-m}} 
\prod_{j=1}^{m}(q^{k-j+1+m}-q^{-(k-j+1+m)}).
\eea

Now, by the symmetry of $F$, we have
\bea
F(-\lambda,\mu)F(\lambda,\nu)=
F(\mu,-\lambda)F(\nu,\lambda).
\eea

Therefore, 
\bea
Res_{\lambda=k}(F(-\lambda,\mu)F(\lambda,\nu))=
Res_{\lambda=k}(F(\nu,\lambda))F(\mu,-k)=
\eea
\bea
F(\nu,-k)\ln(q^{-2})\prod_{j=1, j\neq m+1-k}^{m}
\frac{1}{q^{-k-j+1+m}-q^{k+j-1-m}}
\prod_{j=1}^{m}(q^{k-j+1+m}-q^{-(k-j+1+m)})
F(\mu,-k)=
\eea
\bea
F(\nu,-k)Res_{\lambda=k}(F(\mu,\lambda))=
-F(\nu,-k)Res_{\lambda=-k}(F(\mu,-\lambda))=
-Res_{\lambda=-k}(F(-\lambda,\mu)F(\lambda,\nu)).
\eea

Proposition \ref{resi} and Theorem \ref{chamberindep} for $sl_2$
are proved. 

\subsection{Proof of Theorem \ref{chamberindep} in the dominant case.}

Theorem \ref{chamberindep} claims that for any two big weights
$\xi,\eta\in \hr$, one has $I(\xi)=I(\eta)$. 
In this subsection we will establish the following special case 
of Theorem \ref{chamberindep}: if $\xi$ is a big dominant 
weight, and $s_j$ a simple reflection, then $I(\xi)=I(s_j\xi)$.
We will refer to this case as the dominant case. 

Thus, let us take a big dominant $\xi$. So we have $(\xi,\alpha_j^\vee)=a>>0$.
We may assume, without loss of generality, that 
$(\xi, \alpha_j)<<(\xi,\alpha)$ for
all $\alpha\in \Sigma_+, \alpha \neq \alpha_j$
(this situation may be achieved by moving $\xi$ within its Weyl chamber). 

Let us decompose $\h^*$ as $\Bbb C\alpha_j\oplus \Bbb
(\Bbb C\alpha_j)^\perp$, and write $\lambda\in \h^*$ accordingly as a pair 
$(\lambda_j, \lambda_j^{\perp})$, where
$\lambda_j:=(\lambda,\alpha_j^\vee)\in \Bbb C$. 
Then the cycle $C_\xi$ gets represented as 
a product $L_+\times C_\xi'$, where 
$L_+=(\xi,\alpha_j^\vee)+i\Bbb
R\subset \Bbb C$. Similarly, since 
$s_j\lambda=(-\lambda_j, \lambda_j^{\perp})$,
we have $C_{s_j\xi}=L_-\times C_\xi'$, 
where $L_-=-(\xi,\alpha_j^\vee)+i\Bbb
R$. This means that the integral with respect to $d\lambda$ over
either cycle is representable as an iterated integral. By Fubini's theorem,
 we are free to choose in which order to compute the iterated integral. 
We will choose it so that we integrate first with respect to 
the scalar variable $\lambda_j$. 

Thus, to show that the integral over $C_\xi$ 
is equal to the integral over $C_{s_j\xi}$, it is sufficient to move the
contour of integration with respect to $\lambda_j$ from $L_+$ to
$L_-$ and show that the sum of residues at poles of the integrand
that we pass on the way is zero. Then the equality is guaranteed by the residue
theorem, given that the function $g$ is rapidly decaying.  

So let us study the poles of the integrand. 

\begin{proposition} Let 
$\lambda=(\lambda_j,\lambda_j^\perp)$. Fix $\lambda_j^\perp$ and let 
Let $G(\lambda_j):=F^{V^*}(-\lambda,\mu)^*F^V(\lambda,\nu)$.
Then all the poles of $G(z)$ between the lines 
${\rm Re}(z)=-a$ and ${\rm Re}(z)=a$ are simple,
and located at the points $z=k+l\kappa$, where $l$ is an integer, and 
$k=\pm 1,...,\pm N$. 
\end{proposition} 

\begin{proof}
The proposition follows from the fact that $(\xi,\alpha_j)<<(\xi,\alpha)$
for positive roots $\alpha\ne \alpha_j$, and from Lemma \ref{simppol}. 
\end{proof} 

The function $G$ is $\kappa$-quasi-periodic. 
Therefore, to prove Theorem \ref{chamberindep}, it suffices to establish

\begin{proposition} \label{resi1}
The residue of the function $G(z)$ 
at $z=k$ equals minus its residue at $z=-k$, k=1,...,m.  
\end{proposition} 

The rest of the subsection is occupied by the proof of this
proposition. 

Clearly, it suffices to assume that $\mu,\nu$ are generic. 
Let us restrict the representations $M_\mu$, $M_\nu$ and $V$ 
to the Hopf subalgebra $U_{q^{d_j}}(sl_2)$ of $U_q(\g)$ generated
by $e_j,f_j$ and $q^{\pm h_j}$. Then $M_\mu,M_\nu$ will decompose 
in a direct sum of Verma modules over $U_{q^{d_j}}(sl_2)$, 
and $V$ will decompose into a direct sum of finite
dimensional representations of $U_{q^{d_j}}(sl_2)$. 
Therefore, it is easy to see that the matrix elements of the expression 
$F^{V^*}(-\lambda,\mu)^*F^V(\lambda,\nu)$
are (infinite) linear combinations of 
matrix elements of similar expressions for $U_{q^{d_j}}(sl_2)$, 
with coefficients independent of $z$ 
(i.e. dependent only on $\lambda_i^\perp$).
Hence, the statement in question reduces to the case 
$\g=sl_2$, where it was proved in the previous subsection. 
Proposition \ref{resi1} and Theorem \ref{chamberindep} in the
dominant case are
proved. 

\subsection{Proof of Theorem \ref{chamberindep} in the general case.}

Now we will prove Theorem \ref{chamberindep} in general, using
that it is true in the dominant case. 
This is done using the dynamical Weyl group. Namely, 
by formula (\ref{symmetry}), we have 
\bean\label{weyl1}
F^{V^*}(-\lambda,\mu)^*F^V(\lambda,\nu)=
\eean
\bea
{\mathcal A}_{w,V}(w^{-1}\mu)F^{V^*}(-w^{-1}\lambda,w^{-1}\mu)^*
{\mathcal A}_{w,V^*}(-w^{-1}\lambda)^*
{\mathcal A}_{w,V}(w^{-1}\lambda)
F^{V}(w^{-1}\lambda,w^{-1}\nu){\mathcal A}_{w,V^*}(w^{-1}\nu)^*.
\eea

By Proposition 21 of \cite{EV3},
on $V[0]$ 
one has ${\mathcal A}_{w,V^*}(-w^{-1}\lambda)^*={\mathcal A}_{w^{-1},V}(\lambda)$.
Therefore, the product of two ${\mathcal A}$ operators in the
middle is ${\mathcal A}_{w^{-1},V}(\lambda){\mathcal
A}_{w,V}(w^{-1}\lambda)$, which is
$1$ by the cocycle condition (\ref{cocycle}).
Therefore, formula (\ref{weyl1}) can be rewritten in the form
\bean\label{weyl2}
F^{V^*}(-\lambda,\mu)^*F^V(\lambda,\nu)=
\eean
\bea
{\mathcal A}_{w,V}(w^{-1}\mu)F^{V^*}(-w^{-1}\lambda,w^{-1}\mu)^*
F^{V}(w^{-1}\lambda,w^{-1}\nu){\mathcal A}_{w,V^*}(w^{-1}\nu)^*.
\eea

Now let $\xi$ be a big dominant weight, and $s_j$ a simple
reflection. We have shown in the previous subsection that 
\bea
\int_{C_{\xi}}g(\lambda)F^{V^*}(-\lambda,\mu)^*F^V(\lambda, \nu)d\lambda=
\int_{C_{s_j\xi}}g(\lambda)F^{V^*}(-\lambda,\mu)^*F^V(\lambda,
\nu)d\lambda,
\eea

Therefore, using formula (\ref{weyl2}), we find, $\forall w \in \W$
\bea
{\mathcal A}_{w,V}(w^{-1},\mu)
\int_{C_{\xi}}g(\lambda)F^{V^*}(-w^{-1}\lambda,w^{-1}\mu)^*F^V(w^{-1}\lambda,
w^{-1}\nu)d\lambda\cdot
{\mathcal A}_{w,V^*}(w^{-1}\nu)=\\
{\mathcal A}_{w,V}(w^{-1},\mu)
\int_{C_{s_j\xi}}g(\lambda)F^{V^*}(-w^{-1}\lambda,w^{-1}\mu)^*F^V(w^{-1}\lambda,
w^{-1}\nu)d\lambda\cdot
{\mathcal A}_{w,V^*}(w^{-1}\nu).
\eea
Canceling the ${\mathcal A}$-operators and changing
$w^{-1}\lambda$ to $\lambda$, $w^{-1}\mu$ to $\mu$, and
$w^{-1}\nu$ to $\nu$, we get 
\bea
\int_{C_{w\xi}}g(\lambda)F^{V^*}(-\lambda,\mu)^*F^V(\lambda, \nu)d\lambda=
\int_{C_{ws_j\xi}}g(\lambda)F^{V^*}(-\lambda,\mu)^*F^V(\lambda,
\nu)d\lambda,
\eea
which implies Theorem \ref{chamberindep} in full generality, as
any element of $\W$ can be written as a product of simple
reflections. 

\section{Proof of Theorem \ref{1}}

\subsection{Reduction to big dominant integral $\mu, \nu$}

Let $V$ be a finite dimensional $U_q(\g)$ module.
 From now till the end of this section 
fix $\beta\in P$. Let $\mu\in \h^*$ be generic, and $\nu=\mu-\beta$.

\begin{lemma}\label{rational}
The integral 
\bean
\int _{ C_\xi / {\kappa Q^\vee} }\, F^V(\mu, -\la)\ F^V(\la, \nu)\ d\la\
% =\ \dl_{\mu, \nu}\ Q^{-1}_V (-\mu - \rho) \ .
\notag
\eean
(for big $\xi$) is a rational function of 
the variables $q^{2(\mu, \al_i)}$.
\end{lemma}

\begin{proof}
As explained in Section 2, the matrix elements of the function
$F(\lambda,\mu)$ are representable as finite sums $\sum_{j=1}^m
q^{-2(\lambda,\mu)}
f_j(\lambda)g_j(\mu)$, where $f_j(\nu)$ and $g_j(\nu)$ are
rational functions of the variables $q^{2(\nu,\alpha_i)}$. 
This immediately implies the statement.  
\end{proof}

\begin{corollary}
Theorem \ref{1} is true if it is true for all big anti-dominant
integral $\mu$.
\end{corollary}

\begin{proof}
The corollary follows from the fact that
 a rational function of $q^{2(\mu,\alpha_i)}$ 
is determined by its values at big anti-dominant integral $\mu$. 
\end{proof} 

\subsection{Proof of Theorem \ref{1}}
We continue the proof of Theorem \ref{1}. As was shown in the
previous subsection, we may 
assume that $\mu$ is big anti-dominant integral (then so is
$\nu$, since $\beta$ is fixed and belongs to $P$).

Substituting the generalized Weyl character formula (\ref{dWeyl})
into the left hand side
of the formula of Theorem \ref{1''} we get
\bean\label{nachalo}
LHS=
\eean
\bea
{1\over |\W|}\sum_{w,y\in \W}(-1)^w (-1)^y
\int_{C_\xi/{\kappa Q^\vee}} {\mathcal A}_{y,V}(\mu)^{-1}F^{V^*}(-\la, y\mu)^*
F^{V}(\la, w\nu)({\mathcal A}_{w,V^*}(\nu)^{-1})^* d \la.
\eea
Using the symmetry \Ref{symmetry} and the cocycle condition \Ref{cocycle},
we find
\bean\label{affa=aff}
{\mathcal A}_{y,V}(\mu)^{-1} 
F^{V^*}(-\lambda,y\mu)^*
F^{V}_\nu(\lambda,w\nu)
({\mathcal A}_{w,V^*}(\nu)^{-1})^* 
=
\eean
\bea
{\mathcal A}_{y,V}(\mu)^{-1} 
F^{V^*}(-\lambda,y\mu)^*
{\mathcal A}_{w,V}(w^{-1}\lambda) 
F^{V}(w^{-1}\lambda,\nu)
=
\eea
\bea
{\mathcal A}_{y,V}(\mu)^{-1} 
[{\mathcal A}_{w,V^*}(-w^{-1}\lambda) 
F^{V^*}(-w^{-1}\lambda,w^{-1}y\mu)\times
{\mathcal A}_{w,V}(w^{-1}y\mu)^*]^* 
{\mathcal A}_{w^{-1},V^*}(-\lambda)^* 
F^{V}(w^{-1}\lambda,\nu)
=
\eea
\bea
{\mathcal A}_{y,V}(\mu)^{-1} 
{\mathcal A}_{w,V}(w^{-1}y\mu)
F^{V^*}(-w^{-1}\lambda,w^{-1}y\mu)^* 
F^{V}(w^{-1}\lambda,\nu)
=
\eea
\bea
{\mathcal A}_{w^{-1}y,V}(\mu)^{-1}
F^{V^*}(-w^{-1}\lambda,w^{-1}y\mu)^* 
F^{V}(w^{-1}\lambda,\nu).
\eea
Now we will use Theorem \ref{chamberindep1}. By this theorem and
the last formula, equality (\ref{nachalo})
can be written as
\bean\label{step2}
LHS=\sum_{w\in\ \W} \  (-1)^w\ {\mathcal A}_{w,V}(\mu)^{-1}\
\int_{C_\xi/{\kappa Q^\vee}}\  F^{V^*}(-\la, w\mu)^*\
F^{V}(\la, \nu)\  d \la\ .
\eean
%Now Theorem \ref{1} follows from

Now we have 

\begin{lemma}\label{end-2}
If $w\in\W$ is not the identity element, then
$$
\int_{C_\xi/{\kappa Q^\vee}}\  F^{V^*}(-\la, w\mu)^*\
F^{V}(\la, \nu)\  d\la\ =\ 0 \ .
$$
\end{lemma}
\begin{proof}
It is easy to see that the integral tends to zero
as  $\xi$ tends to infinity inside  the dominant Weyl chamber. 
At the same time the integral does not change its value under
variation of $\xi$. Hence, the integral is zero if $\xi$ is in
the dominant chamber. It follows from Theorem
\ref{chamberindep1} that the same is true for all other Weyl
chambers. The lemma is proved. 
\end{proof}

By Lemma \ref{end-2} and formula \Ref{step2}, 
the left hand side of the formula of Theorem
\ref{1''} is
given by the formula
$$
LHS=\int_{C_\xi/{\kappa Q^\vee}}\  F^{V^*}(-\la, \mu)^*\
F^{V}(\la, \nu)\  d\la\ .
$$
Using the symmetry \Ref{Sym} we see that Theorem \ref{1''}
implies Theorem \ref{1}. Thus Theorem \ref{1} is proved. 

\section{Proof of Theorem \ref{2}}

\subsection{A polynomiality lemma}

\begin{lemma} \label{polyn} 
for generic $\nu \in \h^*$ and big anti-dominant integral $\mu$,
the function \linebreak
$F^{V^*}_\mu(-\lambda)^*F^{V}(\lambda,\nu)$ is a Laurent polynomial of 
$q^{2(\lambda,\alpha_i)}$ multiplied by $q^{2(\lambda,\mu-\nu)},$
\end{lemma}

\begin{proof}
It suffices to check that the same statement is true for 
the function 
$$
(\Psi^v(\lambda,\nu),\Psi_\mu^{v_*}(-\lambda))
\delta_q(-\lambda)\delta_q(\lambda)
$$ 
for big dominant integral
$\mu$. Arguing as in the proof of Theorem \ref{1'},
we find that 
$$
(\Psi^v(\lambda,\nu),\Psi_\mu^{w^*}(-\lambda))=
\Tr[(1\otimes (,)\otimes 1)(\Phi_\nu^v\otimes
\tilde\Phi_\mu^{v_*})
(q^{2\lambda}\otimes q^{2\lambda})],
$$
where the trace is taken in the tensor product $M_\nu\otimes
{}^*L_{\mu}$. This tensor product decomposes as 
$M_\nu\otimes {}^*L_\mu=\oplus H_{\mu\nu^*}^\beta M_\beta$. 
Since the operator $X=(1\otimes (,)\otimes 1)(\Phi_\nu^v\otimes
\tilde\Phi_\mu^{v_*})$ is an intertwiner from 
$M_\nu\otimes {}^*L_\mu$ to itself, as in the proof of Theorem
\ref{1'} we find that
$(\Psi^v(\lambda,\nu),\Psi_\mu^{w^*}(-\lambda))$, as a function of
$\lambda$, is a finite linear combination of characters of Verma
modules $M_\beta$. Hence
$\delta_q(\lambda)(\Psi^v(\lambda,\nu),\Psi_\mu^{w^*}(-\lambda))$ is
equal to $q^{2(\lambda,\mu-\nu)}$ times a Laurent polynomial. This implies
the statement of the lemma.
\end{proof}

\subsection{A modification of Theorem \ref{2}}

\begin{proposition} For a big anti-dominant integral $\mu$, one has
\bean\label{modif}
\frac{1}{|\W|}\int_{C_\xi}F^{V^*}_\mu(-\lambda)^*
\sum_{w\in W}(-1)^w
F^{V}(\lambda,w\nu)
({\mathcal A}_{w,V^*}(\nu)^{-1})^*
q^{-(\lambda,\lambda)}
d\lambda=
\eean
\bea
q^{(\mu,\mu)+(\nu,\nu)}F^{V^*}_\mu(\nu)^*.
\eea
\end{proposition}

\begin{proof}
Let $F_0^V(\lambda,\nu):=F^V(\lambda,\nu)q^{2(\lambda,\nu)}$. 
We have 
\bea
q^{-(\nu,\nu)}LHS=
\frac{1}{|\W|}\int_{C_\xi}F^{V^*}_\mu(-\lambda)^*
\sum_{w\in \W}(-1)^w
F_0^{V}(\lambda,w\nu)
({\mathcal A}_{w,V^*}(\nu)^{-1})^*
q^{-(\lambda+w\nu,\lambda+w\nu)}
d\lambda=\\
\frac{1}{|\W|}\int_{C_{\xi+{\rm Re}(w\nu)}}
\sum_{w\in W}(-1)^w F^{V^*}_\mu(-\lambda+w\nu)^*
F_0^{V}(\lambda-w\nu,w\nu)
({\mathcal A}_{w,V^*}(\nu)^{-1})^*
q^{-(\lambda,\lambda)}
d\lambda.
\eea
By Lemma \ref{polyn}, the cycle of integration 
$C_{\xi+{\rm Re}(w\nu)}$ may be replaced with $C_\xi$,
and furthermore the integral is a rational function of $q^{2(\nu,\alpha_i)}$. 
Therefore, it suffices to check the proposition in the case when
$\nu$ is big anti-dominant integral. But in this case, the result
follows from Theorem \ref{2''} and the generalized Weyl character
formula \Ref{dWeyl}.
\end{proof}

\begin{corollary}\label{expand} 
For big anti-dominant integral $\mu$ and generic $\nu$
\bea
\frac{1}{|\W|}\sum_{w,y\in \W}(-1)^{w}(-1)^{y}\int_{C_\xi}
A_{y,\nu}(\mu)^{-1}F^{V^*}(-\lambda,y\mu)^*F^V(\lambda,w\nu)
(A_{w,V^*}(\nu)^{-1})^*
q^{-(\lambda,\lambda)}d\lambda=
\eea
\bean\label{expand1}
q^{(\mu,\mu)+(\nu,\nu)}\sum_{z\in W}(-1)^z
A_{z,V}(\mu)^{-1}F^{V^*}(\nu,z\mu)^*.
\eean
\end{corollary}

The corollary follows by applying 
the generalized Weyl character formula \Ref{dWeyl} to formula \Ref{modif}. 

\subsection{Proof of Theorem \ref{2}}

Arguing as in the proof of Theorem \ref{1} (using Theorem
\ref{chamberindep} and formula (\ref{affa=aff})), 
we can rewrite identity \Ref{expand1} in the
form
\bean\label{shorter}
q^{-(\mu,\mu)-(\nu,\nu)}\sum_{z\in \W}(-1)^z \int_{C_\xi}
A_{z,V}(\mu)^{-1}F^{V^*}(-\lambda,z\mu)^*
F^V(\lambda,\nu)q^{-(\lambda,\lambda)}d\lambda=
\eean
\bea
\sum_{z \in \W}(-1)^zA_{z,V}(\mu)^{-1}F^{V^*}(\nu,z\mu)^*,
\eea
This identity has been established for big anti-dominant integral
$\mu$. However, we will now use it to prove Theorem \ref{2}, 
which is equivalent to the statement that this identity holds 
for generic complex $\mu$ and, furthermore, term-by-term. 

The proof is based on the following lemma. 
Let $\mu\in \h^*$ be a fixed weight, which does not belong to the pole 
divisor of $F^V(\lambda,\mu)$. 

\begin{lemma} \label{holom}
The function 
\bea
G(\mu,\nu):=\int_{C_\xi} F^{V^*}(-\lambda,\mu)^*F^V(\lambda,\nu)
q^{-(\lambda,\lambda)-(\mu,\mu)-(\nu,\nu)+2(\mu,\nu)}d\lambda.
\eea
is holomorphic in $\nu$ for ${\rm Re}\ \nu$ big anti-dominant, and 
$\kappa Q^\vee$-periodic in $\nu$; that is, $G$ is 
a holomorphic function of $z_j:=q^{-2(\nu,\alpha_j)}$ in the
region $0<|z_j|<q^{dN}$, where $N$ is as in Lemma \ref{simppol}, and
$d$ is the ratio of squared lengths of long and short roots. Furthermore, 
$G$ is meromorphic in the region $|z_j|<q^{dN}$, and 
its Laurent coefficients are rational functions of $q^{2(\mu,\alpha_i)}$.
\end{lemma}

\begin{proof} By Theorem \ref{chamberindep}, we may assume 
that $\xi$ lies in the dominant chamber.
By Lemma \ref{simppol}, the trace function 
$F^V(\lambda,\mu)$ may be written as 
\bean\label{ratio}
F^V(\lambda,\mu)=q^{-2(\lambda,\mu)}\frac{\tilde{F}^V(\lambda,\mu)}{\prod_{\alpha\in
\Sigma_+}
\prod_{k=1}^N(1-q^{2(\lambda,\alpha)-k(\alpha,\alpha)})},
\eean
where $\tilde{F}^V(\lambda,\mu)$
is a trigonometric polynomial in $\lambda$ whose coefficients
are rational functions of $q^{2(\mu,\alpha_i)}.$

Therefore, we have
\bea
G(\mu,\nu)=\int_{C_\xi}
q^{-(\lambda-\mu+\nu,\lambda-\mu+\nu)}
\frac{{H}(\lambda,\mu,\nu)}
{\prod_{\alpha\in \Sigma_+}\prod_{k=1}^N(1-q^{2(\alpha,\lambda)-k(\alpha,\alpha)})
(1-q^{2(\alpha,\lambda)+k(\alpha,\alpha)})}d\lambda,
\eea
where $H$ is a trigonometric polynomial in $\lambda$ 
whose coefficients are rational in $q^{2(\mu,\alpha_j)}$ and $q^{2(\nu,\alpha_j)}$

Now make a change of variable 
$\eta=\lambda-\mu+\nu$. Then we get
\bea
G(\mu,\nu)=\int_{C_{\xi-{\rm Re}(\mu-\nu)}}
q^{-(\eta,\eta)}
\frac{H(\eta+\mu-\nu,\mu,\nu)}
{\prod_{\alpha\in \Sigma_+}\prod_{k=1}^N(1-q^{2(\alpha,\eta+\mu-\nu)-k(\alpha,\alpha)})
(1-q^{2(\alpha,\eta+\mu-\nu)+k(\alpha,\alpha)})}d\eta.
\eea
The fractions $\frac{1}{1-q^{2(\alpha,\eta+\mu-\nu)\pm k(\alpha,\alpha)}}$ 
 in the integrand can be expanded into geometric series 
for ${\rm Re}(\alpha_i,\eta+\mu-\nu)>>0$:
$$
\frac{1}{1-q^{2(\alpha,\eta+\mu-\nu)\pm k(\alpha,\alpha)}}=
\sum_{m\ge 0}q^{m(2(\alpha,\eta+\mu-\nu)\pm k(\alpha,\alpha))}.
$$
 The inequalities ${\rm Re}(\alpha_i,\eta+\mu-\nu)>>0$ are satisfied 
on the cycle of integration, as $\xi$ is a big dominant weight.
Thus we can expand the integrand in a series and then 
integrate the series term-wise. This shows that we can expand 
$G(\mu,\nu)$ in a Laurent series with respect
to $z_j=q^{-2(\nu,\alpha_j)}$, which is convergent
for small enough $z_j$, i.e. 
for big antidominant ${\rm Re}\ \nu$
(in fact, it is easy to see that $|z_j|<q^{dN}$ is 
sufficient for convergence). 
The coefficients of this series are 
integrals of trigonometric polynomials 
of $\eta$ with coefficients rational in $q^{2(\mu,\alpha_i)}$ 
against the measure $q^{-(\eta.\eta)}d\eta$. 
Thus they are rational functions in $q^{2(\mu,\alpha_i)}$, and
the lemma is proved. 
\end{proof}

Now we can finish the proof of Theorem \ref{2}.
Lemma \ref{holom} implies that both sides of identity
\Ref{shorter} admit Laurent expansions with respect to the
variables $z_j=q^{-2(\nu,\alpha_j)}$ for $|z_j|<q^{dN}$. 
Let $\beta$ be a fixed integral weight, and 
$C_{\beta,l}(\mu),C_{\beta,r}(\mu)$ be the coefficients
of $q^{-2(\nu,\beta)}$ in this expansion on the left,
respectively right hand side of identity \Ref{shorter}. 
Let also $C_{\beta,l}^1(\mu),C_{\beta,r}^1(\mu)$
be the coefficients of the same term in the series on the left
and right hand sides of \Ref{shorter} given exclusively by the terms with 
the Weyl group element $z$ equal to $1$.
It is easy to see that if $\mu$ is a big anti-dominant integral
weight, then $C_{\beta,l}(\mu)=C_{\beta,l}^1(\mu)$, 
$C_{\beta,r}(\mu)=C_{\beta,r}^1(\mu)$
(i.e., terms with $z\ne 1$ contribute only to very high terms of the
expansion). On the other hand, since identity 
\Ref{shorter} is known to hold for big anti-dominant integral
$\mu$, we have for such $\mu$:
$C_{\beta,l}(\mu)=C_{\beta,r}(\mu)$. 
Therefore, we conclude that for big anti-dominant integral $\mu$,
one has $C_{\beta,l}^1(\mu)=C_{\beta,r}^1(\mu)$. 
But by Lemma \ref{holom}, both $C_{\beta,l}^1(\mu)$ and
$C_{\beta,r}^1(\mu)$
are rational functions of $q^{2(\mu,\alpha_i)}$. 
Therefore, the equality $C_{\beta,l}^1(\mu)=C_{\beta,r}^1(\mu)$
holds for generic $\mu$, and hence the terms 
on the left and right sides of \Ref{shorter} 
corresponding to $z=1$ are equal: 
\bea
q^{-(\mu,\mu)-(\nu,\nu)}\int_{C_\xi}
F^{V^*}(-\lambda,\mu)^*
F^V(\lambda,\nu)q^{-(\lambda,\lambda)}d\lambda=
F^{V^*}(\nu,\mu)^*.
\eea
This implies Theorem \ref{2}, as $F^{V^*}(\nu,\mu)^*=F^V(\mu,\nu)$.

\section{Integral Transforms}

In this section we prove Theorem \ref{transform}.

\subsection{Proof that the integral transforms are well defined
and continuous}
For any  vectors $v\in V$ and $v_* \in V^*$, we have
$$
(F^V(\la, \mu)v, v_*)\ = \ q^{-2(\la, \mu)}\ \sum_j\ f_j(\la)\ g_j(\mu)
$$
for suitable functions $f_j(\la), g_j(\la)$. Each of these functions
has the form
\bean\label{newton}
{P(\la)\over \prod_{\al\in \Sigma_+} 
\prod_{k=1}^N (1-q^{2(\la,\al)-k(\alpha,\alpha)})}\ ,
\eean
where $P(\la)$ is a Laurent polynomial in variables $q^{2(\al_i, \la)}$,
whose Newton polyhedron is contained in the Newton 
polyhedron of the denominator (this follows from \cite{ESt},
Proposition 2.2).

For any  vectors $v\in V$ and $v_* \in V^*$, the function 
$(Q (-\la-\rho) v, v_*)$ also has such a form (again, from
\cite{ESt}).

\begin{lemma}\label{bounded}
Let a function $f(\la)$ be of the form \Ref{newton}. Let $\xi\in\hr^*$ be big.
Then all derivatives of the restriction of the 
function $f(\la)$ to ${C_\xi}$ or $D_\eta$ are
bounded from above.
\end{lemma}

\begin{proof}
The statement about the restriction to $C_\xi$ is obvious, since in
this case the function is periodic with respect to a lattice. 
So we need to prove only the statement about restriction to 
$D_\eta$. For this purpose, take $\lambda=i\eta+y$, $y\in
\h_{\Bbb R}$. Then the binomial factor
$(1-q^{2(\la,\al)-k(\alpha,\alpha)})$
from formula \Ref{newton} satisfies the following lower bound: 
$$
|1-q^{2(\la,\al)-k(\alpha,\alpha)}|\ge C(1+q^{2(y,\alpha)-k(\alpha,\alpha)})
$$
for some $C=C(\eta)>0$
(This follows from the elementary estimate 
$$
|1-ae^{i\theta}|\ge
(1+a)|\sin(\theta/2)|
$$ 
if $a>0$). This lower bound implies the
claim, since it implies that the denominator of \Ref{newton}
is bounded from below by a polynomial with the same Newton
polyhedron and positive coefficients. 
\end{proof}

Lemma \ref{bounded} implies that the integral transforms of Theorem \ref{transform}
are well defined and continuous. Indeed, let for example $a\in
{\mathcal S}_\eta(C_\xi)\otimes V[0]$, $a(\lambda)=\sum a_m(\lambda)v_m$, where 
$v_m$ is a basis of $V[0]$. Then 
the components of the vector function $K_{\rm Im}^Va$
are finite linear combinations of functions of the form 
\bea
b(\mu)=g(\mu)\int_{C_\xi}q^{2(\lambda,\mu)}f(\lambda)a(\lambda)d\lambda
\eea
where $f$ and $g$ are of the form \Ref{newton}.
By Lemma \ref{bounded}, the multiplication operators by $f$ and
$g$ preserve the spaces of Schwartz functions and are continuous. So the fact that 
$K_{\rm Im}^V$ is well defined and continuous follows from the classical
fact that the Fourier transform is a continuous isomorphism 
of the Schwartz space to itself. The same proof applies to $K_{\rm Re}^V$.

\subsection{Proof that $K^V_{\IM} \ K^V_{\RE} \ = \text{Id}$.}

\begin{proposition}\label{imre}
We have
$$
K^V_{\IM} \ K^V_{\RE} \ = \text{Id}\ .
$$
\end{proposition}

\begin{proof} We will use the following standard lemma 
from Fourier analysis. 

Let $f,g$ be smooth functions on $\Bbb R^n\times \Bbb R^n$ 
such that $f(x,y)=e^{i(x,y)}f_0(x,y)$, 
$g(x,y)=e^{i(x,y)}g_0(x,y)$, 
and $f_0,g_0$ are periodic with respect to the lattice $(2\pi\Bbb
Z)^n$ in the first variable. 
Assume that for any $y,z\in \Bbb R^n$ 
such that $y-z\in (2\pi\Bbb Z)^n$, 
one has
$$
\int_{(\Bbb R/2\pi\Bbb Z)^n}f(x,y)g(-x,z)dx=\delta_{0,y-z}, 
$$
Define integral 
transforms $K_1: C^\infty_0(\Bbb R^n)\to C^\infty(\Bbb
R^n)$, $K_2: {\mathcal S}(\Bbb R^n)\to C^\infty(\Bbb R^n)$ by  
$$
(K_1a)(x)=\int_{\Bbb R^n}g(-x,z)a(z)dz,
\qquad
(K_2b)(y)=\int_{\Bbb R^n}f(x,y)b(x)dx
$$
(here $C^\infty_0$ is the space of functions with compact support).

\begin{lemma}\label{inve}
Under these conditions, $K_1$ takes values in ${\mathcal S}(\Bbb
R^n)$, and $K_2\circ K_1={\rm Id}$
(for an appropriate normalization of $dx$ and $dz$). 
\end{lemma} 

The proof of the lemma is standard, using Fourier series
expansions, and will be omitted.

Lemma \ref{inve} and Theorem \ref{1} imply Proposition
\ref{imre}, since matrix elements of $F^V(\lambda,\mu)$, 
$\lambda\in C_\xi$, $\mu\in D_\eta$, are of the same type as the
functions $f(x,y)$, $g(x,y)$ considered above. 
\end{proof}

\subsection{Proof that $K^V_{\RE} \ K^V_{\IM} \ = \text{Id}$.}

\begin{lemma}\label{vs}
Let $X, Y$ be vector spaces and $A : X \to Y$ and $B : Y \to X$ linear operators.
Assume that the operator $A$ is injective, and $A B \ =\ \text{Id}$.
Then $B A \ = \ \text{Id}$.
\end{lemma}

\begin{proof} For every $y\in Y$ we have $A ( B A y - y ) \ = \
0$. Since $A$ is injective,
we have $B A \ =\ \text{Id}$.
\end{proof}

Lemma \ref{vs} shows that in order to finish the proof of Theorem
\ref{transform} 
it is enough to prove
the following 

\begin{proposition}\label{last}
Under the conditions of Theorem \ref{transform} the operator
$K^V_{\IM}$ is injective.
\end{proposition}

\begin{proof} 
First of all, using the dynamical Weyl group symmetry, we can
reduce the problem to the case when $\xi$ is anti-dominant. 

Suppose that $a\in {\mathcal S}_\eta(C_\xi)$, and 
$K^V_{\IM}a=\int_{C_\xi}F^V(\mu,-\la)a(\la)d\la=0$. 
As we know (see e.g. formula \Ref{ratio}), 
the function $F^V$ is representable in the form 
$F^V(\mu,\nu)=q^{-2(\mu,\nu)}\frac{\hat F^V(\mu,\nu)}{Z(\nu)}$, where 
$Z(\nu)$ is a product of binomial terms, and $\hat F$ is a
trigonometric polynomial in $\nu$. 
This means, 
\bean\label{hateqn}
\int_{C_\xi}\hat q^{-2(\mu,\lambda)}\hat F^V(\mu,-\lambda)b(\lambda)=0,
\eean
where $b(\lambda)=Z(-\lambda)a(\lambda)\in {\mathcal S}_\eta(C_\xi)$
(as $Z$ is periodic). 

The function 
$\hat F^V(\mu,-\lambda)$ can be written in the form 
\bean\label{1plus}
\hat F^V(\mu,-\lambda)=1+\sum_{\beta\in L}E_\beta(\mu)q^{2(\beta,\lambda)},
\eean
where $L\subset Q_+\setminus \lbrace{0\rbrace}$ 
is a finite set, and $E_\beta(\mu)$ are bounded
rational functions of $q^{2(\mu,\alpha_i)}$ on $D_\eta$.

Using equation \Ref{1plus}, 
equation \Ref{hateqn} can be written
in the form of a difference equation:
\bean\label{hateqn1}
\hat b(\mu)=-\sum_{\beta\in L}E_\beta(\mu)\hat b(\mu+\beta).
\eean
If $\xi$ is sufficiently big, $\hat b(\mu+\beta)$,
decays rapidly with $\beta\to \infty$
in $Q_+$ (since $\hat b\in S_\xi(D_\eta)$). Thus, solving  
equation \Ref{hateqn} recursively, and using that $E_\beta(\mu)$ 
are bounded, we find that for sufficiently big $\xi$, 
any solution $\hat b\in S_\xi(D_\eta)$ of 
the difference equation \Ref{hateqn1} must be zero. 
The proposition is proved.
\end{proof}

\section{Proofs of Theorems \ref{selfad} and \ref{realint}}

\subsection{Proof of Theorem \ref{selfad}}
Recall that the trace functions are eigenfunctions of ${\mathcal
D}_{U,V}$, namely
\bea
{\mathcal D}_{U,V}^{(\lambda)}F(\lambda,\nu)=\chi_U(q^{-2\nu})F(\lambda,\nu).
\eea
Now let $\mu - \nu \in P.$ Then by Theorem \ref{1},
\bea
\int_{C_\xi/\kappa Q^\vee}
F(\mu,-\lambda)[{\mathcal D}_{U,V}^{(\lambda)}-({\mathcal D}_{U,V^*}^{(\lambda)})^*]
F(\lambda,\nu)d\lambda=
\eea
\bea
\chi_U(q^{-2\nu})\int_{C_\xi/\kappa Q^\vee}F(\mu,-\lambda)
F(\lambda,\nu)d\lambda-
\chi_U(q^{-2\mu})\int_{C_\xi/\kappa Q^\vee}F(\mu,-\lambda)
F(\lambda,\nu)d\lambda=
\eea
\bea
(\chi_U(q^{-2\nu})-
\chi_U(q^{-2\mu}))
\delta_{\mu\nu}Q^{-1}(-\mu-\rho)=0,
\eea
because either $\delta_{\mu\nu}=0$ or $(\chi_U(q^{-2\nu})-
\chi_U(q^{-2\mu}))=0.$

By Lemma \ref{last} this implies that
\bea
({\mathcal D}_{U,V}^{(\lambda)}-({\mathcal D}_{U,V}^{(\lambda)})^*)F(\lambda,\nu)=0.
\eea
This means that 
\bea
\int_{C_\xi/\kappa Q^\vee}q^{-2(\lambda,\mu)}
({\mathcal D}_{U,V}^{(\lambda)}-({\mathcal D}_{U,V}^{(\lambda)})^*)
F(\lambda,\nu)=0.
\eea
By Lemma \ref{last} this means that
\bea
({\mathcal D}_{U,V}-{\mathcal D}_{U,V}^*)q^{2(\lambda,\nu)}=0,
\eea
which easily implies
\bea
{\mathcal D}_{U,V}={\mathcal D}_{U,V}^*.
\eea
The theorem is proved. 

\subsection{Proof of Theorem \ref{realint}}

By Theorem \ref{2} we have
\bea
\int_{C_\xi}F^V(\mu,-\lambda)
F^V(\lambda,\nu)q^{-(\lambda,\lambda)}d\lambda=
q^{(\mu,\mu)+(\nu,\nu)}
F^V(\mu,\nu).
\eea
This is equivalent to saying that for a fixed $\nu$, one has
\bea
K_{\rm Im}^V(q^{-(\lambda,\lambda)}F^V(\lambda,\nu))=
q^{(\nu,\nu)}(q^{(\mu,\mu)}F^V(\mu,\nu))
\eea
Applying $K_{\rm Re}$ to both sides, and 
using Theorem \ref{transform}, we obtain
\bea
q^{-(\lambda,\lambda)}F^V(\lambda,\nu)=
\int_{D_\eta}F^V(\lambda,\mu)Q(-\mu-\rho)q^{(\nu,\nu)+(\mu,\mu)}
F^V(\mu,\nu),
\eea
as desired.


\begin{thebibliography}
\normalsize


\bibitem[Cha]{Cha}
O. Chalykh, 
Macdonald polynomials and algebraic integrability. 
Adv. Math. 166 (2002), no. 2, 193--259.

\bibitem [Ch1]{Ch1} I. Cherednik, 
Difference Macdonald-Mehta conjecture. 
Internat. Math. Res. Notices 1997, no. 10, 449--467.


\bibitem [Ch2]{Ch2} I. Cherednik, 
Macdonald's 
evaluation conjectures and difference 
Fourier transform. Invent. Math. 122 (1995), no. 1, 119--145.

\bibitem [Dr]{Dr} 
V. G. Drinfeld, 
On almost cocommutative Hopf algebras, Leningrad Math. J. 1 (1990), 321--342.


\bibitem [EKi1]{EKi1}
P. Etingof; A. Kirillov, Jr. 
On Cherednik-Macdonald-Mehta identities. 
Electron. Res. Announc. Amer. Math. Soc. 4 (1998), 43--47

\bibitem [EKi2]{EKi2}
P. Etingof and A. Kirillov, Jr., Macdonald's polynomials and
representations of quantum groups, Math. Res. Let. 1 (1994), 279-296.

\bibitem [EKi3]{EKi3}
P. Etingof and A. Kirillov, Jr.,
Representation-theoretic proof of the inner product and
symmetry identities for Macdonald's polynomials,
Compos. Math. 102, (1996), 179-202.

\bibitem[ESt]{ESt} P.Etingof, K.Styrkas,
{\it Algebraic integrability of Macdonald operators and
representations of quantum groups},
Comp. Math., v. 114, p.125-152, 1998.

\bibitem [ESV]{ESV}, P. Etingof, O. Schiffmann, A. Varchenko,  
Traces of intertwiners for quantum groups and difference
equations, II, math.QA/0207157, to appear in LMP.  

\bibitem[EV1]{EV1} P.Etingof, A.Varchenko,
{\it Exchange Dynamical Quantum Groups}, QA/9801135,
Commun.\ Math.\ Phys.\ {\bf 205} (1999), 19--52.

\bibitem[EV2]{EV2} P.Etingof, A.Varchenko,
{\it Traces of Intertwiners for Quantum Groups and Difference Equations, I},
QA/9907181, Duke \ Math.\ Journal {bf 104} (2000), No. 3, 391--432.

\bibitem[EV3]{EV3} P.Etingof, A.Varchenko,
{\it Dynamical Weyl Groups and Applications},
QA/0011001, Advances in Mathematics, v.167, p.74-127 (2002).

\bibitem [FTV1]{FTV1} 
 G.Felder, V.Tarasov, A.Varchenko,
Solutions of the elliptic QKZB equations
   and Bethe ansatz I, q-alg/9606005, in: 
Topics in Singularity Theory, V.I.Arnold's 60th Anniversary
Collection, Advances in the Mathematical Sciences -34, 
AMS Translations, Series 2, v. 180, pp.
45-75, 1997.

\bibitem [FTV2]{FTV2} G. Felder, V. Tarasov, A. Varchenko,
 Monodromy of solutions of the elliptic quantum 
Knizhnik-Zamolodchikov-Bernard difference equations. Internat. J. Math. 10 (1999),
  no. 8, 943--975.

\bibitem [FV1]{FV1} G. Felder, A. Varchenko,
The $q$-deformed Knizhnik-Zamolodchikov-Bernard heat equation. 
Comm. Math. Phys. 221 (2001), no. 3, 549--571.

\bibitem [FV2]{FV2}
G. Felder, A. Varchenko, 
$q$-deformed KZB heat equation: 
completeness, modular properties and ${\rm SL}(3,\Bbb Z)$. Adv. Math. 171 (2002), no. 2,
  228--275.

\bibitem [FV3]{FV3}
G. Felder, A. Varchenko, The elliptic gamma
function and ${\rm SL}(3, Z)\ltimes Z\sp 3$. Adv. Math. 156
(2000), no. 1, 44--76.

\bibitem [Kos]{Kos} 
B. Kostant, On Macdonald's $\eta$-function formula, 
the Laplacian and generalized exponents, Advances in Math. 20 (1976),
179-212. 

\bibitem [M]{M}
I. G. Macdonald, A new class of symmetric functions, Sem. Lothar. Combin.
20 (1988).

\bibitem[STV]{STV} K. Styrkas, V. Tarasov, and A. Varchenko,
{\it How to regularize singular vectors and kill the dynamical Weyl group},
math.QA/0206294.

\bibitem[TV]{TV} V.Tarasov, A.Varchenko, {\it Difference Equations Compatible 
with Trigonometric KZ Differential Equations}, 
QA/0002132,  IMRN 2000, No. 15, 801-829.



\end{thebibliography}
\end{document}